\documentclass[a4paper,10pt,12pt]{article}
\usepackage{amsmath, amssymb, latexsym, amscd, amsthm,amsfonts,amstext}
\usepackage[mathscr]{eucal}
\usepackage{graphicx}
\usepackage{subfig}
\usepackage{float}
\usepackage{color}
\usepackage{hyperref}

\newtheorem{lemma}{Lemma}[section]
\newtheorem{problem}{Problem}[section]

\newtheorem{algorithm}{Algorithm}[section]

\numberwithin{equation}{section}
\input{tcilatex}

\newcommand{\width}{5cm}
\newcommand{\height}{4.5cm}

\begin{document}

\title{A globally convergent method for a 3-D inverse medium problem for the
generalized Helmholtz equation}
\author{Michael V. Klibanov\thanks{%
corresponding author}, Hui Liu and Loc H. Nguyen \\
\\
Department of Mathematics \& Statistics \\
University of North Carolina at Charlotte\\
Charlotte, NC 28223, USA \\
Emails: mklibanv@uncc.edu, hliu34@uncc.edu, lnguye50@uncc.edu}
\date{}
\maketitle

\begin{abstract}
A 3-D inverse medium problem in the frequency domain is considered. Another
name for this problem is Coefficient Inverse Problem. The goal is to
reconstruct spatially distributed dielectric constants from scattering data.
Potential applications are in detection and identification of explosive-like
targets. A single incident plane wave and multiple frequencies are used. A
new numerical method is proposed. A theorem is proved, which claims that a
small neigborhood of the exact solution of that problem is reached by this
method without any advanced knowledge of that neighborhood. We call this
property of that numerical method \textquotedblleft global convergence".
Results of numerical experiments for the case of the backscattering data are
presented. 
\end{abstract}

\textbf{Key Words}: global convergence, coefficient inverse problem,
frequency domain

\textbf{2010 Mathematics Subject Classification:} 35R30.

\section{Introduction}

\label{sec:1}

Potential applications of the Inverse Medium Problem of this paper are in
detection and identification of explosive-like targets using measurements of
electromagnetic data. In the case of time dependent experimental data, this
application was addressed in \cite{BK1,KFB,KSNF1,TBKF1,TBKF2}. In the
current paper, so as in \cite{BK1,KFB,KSNF1,TBKF1,TBKF2}, we calculate
dielectric constants of targets for the case of the frequency dependent
data. Of course, estimates of dielectric constants alone cannot
differentiate between explosives and the clutter. On the other hand, the
radar community is relying now only on the intensity of radar images \cite%
{KSNF1}. Thus, we hope that the additional information about dielectric
constants might lead in the future to the development of algorithms, which
would better differentiate between explosives and the clutter.

An Inverse Medium Problem is the problem of determining one of coefficients
of a PDE from boundary measurements. Another name for it is Coefficient
Inverse Problem (CIP) or Inverse Scattering Problem. We are interested in a
CIP for a generalized Helmholtz equation with the data resulting from a
single measurement event. In other words, the boundary data are generated by
a single direction of the incident plane wave and boundary measurements are
conducted on an interval of frequencies. Thus, we use the minimal number of
measurements for a CIP in the frequency domain. We call a numerical method
for a CIP \emph{globally convergent} if a theorem is proved, which claims
that this method delivers at least one point in a sufficiently small
neighborhood of the exact solution without any advanced knowledge of this
neighborhood.

Currently there exist two types of globally convergent numerical methods for
CIPs with single measurement data. The method of the first type was
completely verified on electromagnetic experimental data, see, e.g. \cite%
{BK1,KFB,KSNF1,TBKF1,TBKF2}. As to the method of the second type, it was
initiated in \cite{Klib97} with a recently renewed interest \cite{BK2,KNT,KK}%
. In particular, in \cite{KNT} numerical experiments are presented.

Methods of both types start from a CIP for a hyperbolic PDE. Next, the
Laplace transform is applied with respect to time. It transforms the
original hyperbolic PDE in the equation 
\begin{equation}
\Delta w-s^{2}c\left( x\right) w=-\delta \left( x-x_{0}\right) ,x\in \mathbb{%
R}^{3},s>0  \label{1.1}
\end{equation}%
with the unknown coefficient $c\left( x\right) >0.$ Using the maximum
principle, one can prove that $w>0$. Next, the function $\partial _{s}\left(
\ln w/s^{2}\right) $ is considered and an integral differential PDE is
obtained for this function. Integration is carried out from $s$ to $\infty .$
In the method of the first type, one truncates those integrals at a
sufficiently large value $s:=\overline{s}.$ Next, one obtains a sequence of
Dirichlet boundary value problems for elliptic PDEs. Solving those PDEs
sequentially as well as updating residuals of those truncated integrals, one
obtains points in a sufficiently small neighborhood of the exact solution
without any advanced knowledge of that neighborhood. This corresponds to the
above definition of the global convergence.

In the method of the second type, one does not truncate that integral.
Rather, one uses Laugerre functions as well as Carleman Weight Function to
construct a Tikhonov-like cost functional, which is strictly convex on any
reasonable bounded set in a Sobolev space. This ensures the convergence of
the gradient method to the unique minimum of that functional starting from
any point of that bounded set. Convergence of minimizers to the exact
solution when the level of the error in the data tends to zero is also
guaranteed.

In this paper we develop a frequency domain analog of the globally
convergent numerical method of the first type. The reason of this is that
one can choose either of two types of measurements for the above application
to detection and identification of explosives: either measurements of time
dependent data, as in \cite{BK1,KFB,TBKF1,TBKF2}, or measurements of
frequency dependent data for a certain interval of frequencies.

One of the most difficult questions to address in this paper is that we need
to work now with the complex valued analog of the function $w>0$ in (\ref%
{1.1}). Let $\widetilde{w}$ be that analog. It is not immediately clear how
to define $\func{Im}\left( \log \widetilde{w}\right) .$ To handle this
difficulty, we modify our previous algorithm of \cite{BK1,KFB,
KSNF1,TBKF1,TBKF2}, using the fact that $\partial \log \widetilde{w}%
=\partial \widetilde{w}/\widetilde{w}.$ So, we use only derivatives of $\log 
\widetilde{w}$. Moreover, the use of those derivatives leads us to a new
scheme of the numerical method, as compared with the one of \cite%
{BK1,KFB,KSNF1,TBKF1,TBKF2}. The second difficult question to address here,
as compared with \cite{BK1,KFB,TBKF1,TBKF2}, is that, unlike (\ref{1.1}),
the maximum principle does not work for the generalized Helmholtz equation.

Globally convergent numerical methods for CIPs for the case of the data
resulting from multiple measurements were developed in \cite{Kab1,Kab2}. We
also refer to the survey \cite{Bao} for numerical methods for CIPs in the
frequency domain with multiple frequencies and to, e.g. \cite%
{Li1,Li2,Li3,Nov1,Nov2} for some other inverse scattering problems in the
frequency domain.

In Section \ref{sec:2} we formulate forward and inverse problems which we
consider. In Section \ref{sec:3} we consider the asymptotic behavior of the
solution of the forward problem when the frequency tends to infinity. In
Section \ref{sec:4} we use the Lippmann-Schwinger equation to establish some
properties of the solution of the forward problem. In Section \ref{sec:5} we
describe our numerical method. In Section \ref{sec:6} we establish existence
and uniqueness theorem of a certain auxiliary boundary value problem.
Section \ref{sec:7} is devoted to the convergence analysis. Numerical
implementation of our method and numerical experiments are described in
Section \ref{sec:8}. We briefly summarize our results in Section 9.

\section{The statement of the inverse scattering problem in the frequency
domain}

\label{sec:2}

Let $B\left( R\right) =\left\{ \left\vert x\right\vert <R\right\} \subset 
\mathbb{R}^{3}$ be the ball of the radius $R$ centered at $0$. Let $\Omega
_{1}\Subset \Omega \Subset B\left( R\right) $ be two domains with boundaries 
$\partial \Omega $ and $\partial \Omega _{1}$ and let the domain $\Omega $
be convex. Both boundaries belong to the class $C^{2+\alpha }$ for some $%
\alpha \in (0,1).$ Here and below $C^{m+\alpha }$ are H\"{o}lder spaces of
complex valued functions, where $m\geq 0$ is an integer. For any domain $%
Q\subseteq \mathbb{R}^{3}$ with the $C^{m+\alpha }$ boundary $\partial Q$
the norm in $C^{m+\alpha }\left( \overline{Q}\right) $ of a complex valued
function $v$ is defined in the natural manner as $\left\Vert v\right\Vert
_{C^{m+\alpha }\left( \overline{Q}\right) }=\left\Vert \func{Re}v\right\Vert
_{C^{m+\alpha }\left( \overline{Q}\right) }+\left\Vert \func{Im}v\right\Vert
_{C^{m+\alpha }\left( \overline{Q}\right) }.$ If a function $v\in
C^{m+\alpha }\left( \mathbb{R}^{3}\right) ,$ then we denote $\left\Vert
v\right\Vert _{m+\alpha }=\left\Vert v\right\Vert _{C^{m+\alpha }\left( 
\mathbb{R}^{3}\right) }.$ We denote norms in the spaces $C^{m+\alpha }\left( 
\overline{\Omega }\right) $ as $\left\vert v\right\vert _{m+\alpha },\forall
v\in C^{m+\alpha }\left( \overline{\Omega }\right) .$ If two functions $%
f,g\in C^{\alpha }\left( \overline{\Omega }\right) ,$ then obviously $%
\left\vert fg\right\vert _{\alpha }\leq \left\vert f\right\vert _{\alpha
}\left\vert g\right\vert _{\alpha }.$ For any complex valued function $f\in
C^{m+\alpha }\left( \overline{\Omega }\right) ,$ we define 
\begin{equation*}
\left\vert \nabla f\right\vert _{m+\alpha }=\dsum\limits_{j=1}^{3}\left\vert
f_{x_{j}}\right\vert _{m+\alpha }.
\end{equation*}

Assume that the spatially distributed dielectric constant $c(x),x\in \mathbb{%
R}^{3}$ satisfies the following conditions: 
\begin{equation}
c(x)\in C^{15}(\mathbb{R}^{3}),\quad c(x)=1+\beta (x),  \label{2.0}
\end{equation}%
\begin{equation}
\beta (x)\geq 0,\text{ }\beta (x)=0\quad \text{for }\>x\in \mathbb{R}%
^{3}\setminus \Omega _{1}.  \label{2.2}
\end{equation}%
The $C^{15}-$smoothness of the function $c\left( x\right) $ was used in \cite%
{KR} for the proof of an analog of Lemma 3.1 (Section \ref{sec:3}). We
consider the following generalized Helmholtz equation 
\begin{equation}
\Delta u+k^{2}c(x)u=0,\quad x\in \mathbb{R}^{3},  \label{2.3}
\end{equation}%
where $u\left( x,k\right) $ is the complex valued wave field and $k>0$ is
the frequency. Let the incident plane wave $u_{0}\left( x,k\right) =\exp
\left( -ikx_{3}\right) $ propagates along the positive direction of the $%
x_{3}-$axis. The total wave field 
\begin{equation}
u\left( x,k\right) =u_{0}\left( x,k\right) +u_{\mathrm{sc}}\left( x,k\right)
,  \label{2.4}
\end{equation}%
is the solution of equation (\ref{2.3}), which satisfies the radiation
condition at the infinity, 
\begin{equation}
\frac{\partial u_{\mathrm{sc}}}{\partial r}+iku_{\mathrm{sc}%
}=o(r^{-1}),\>\>r=|x|\rightarrow \infty .  \label{2.5}
\end{equation}%
Here $u_{\mathrm{sc}}(x,k)$ denotes the scattering wave. It is well known
that the problem (\ref{2.3})-(\ref{2.5}) has unique solution $u\left(
x,k\right) \in C^{2+\alpha }\left( \mathbb{R}^{3}\right) ,$ see Theorems 8.3
and 8.7 in \cite{CK} as well as Section 4.\ Furthermore, Theorem 6.17 of 
\cite{GT} implies that the function $u\left( x,k\right) \in C^{16+\alpha
}\left( \mathbb{R}^{3}\right) .$ We consider the following inverse problem: 
%

\begin{problem}[Coefficient Inverse Problem (CIP)]
\label{Problem 2.1} Let $\underline{k}$ and $\overline{k}$ be two constants
such that $0<\underline{k}<\overline{k}.$ Assume that the function $g\left(
x,k\right) $ is known, where%
\begin{equation}
g\left( x,k\right) =u\left( x,k\right) ,\text{ \ }x\in \partial \Omega ,k\in %
\left[ \underline{k},\overline{k}\right] .  \label{2.6}
\end{equation}%
Determine the function $\beta \left( x\right) $ for $x\in \Omega .$
\end{problem}

Since $c(x)=1$ in $\mathbb{R}^{3}\setminus \Omega $, the function $%
u_{sc}(x,k)$ solves the following problem outside of the domain $\Omega $ 
\begin{equation}
\begin{array}{rcll}
\Delta u_{\mathrm{sc}}+k^{2}u_{\mathrm{sc}} & = & 0 & \mbox{in }\mathbb{R}%
^{3}\setminus \overline{\Omega }, \\ 
u_{\mathrm{sc}} & = & g-u_{0} & \mbox{on }\partial \Omega , \\ 
\partial _{r}u_{\mathrm{sc}}+iku_{\mathrm{sc}} & = & o(r^{-1}) & \mbox{as }%
r\rightarrow \infty .%
\end{array}
\label{2.7}
\end{equation}%
The problem (\ref{2.7}) has unique solution $u_{\mathrm{sc}}\in C^{2}\left( 
\mathbb{R}^{3}\setminus \overline{\Omega }\right) \cap C\left( \mathbb{R}%
^{3}\setminus \Omega \right) $, see Lemma 3.8 and Theorem 3.9 in \cite{CK}.
Also, since we have established above that $u\left( x,k\right) \in
C^{2+\alpha }\left( \mathbb{R}^{3}\right) ,$ then $u_{\mathrm{sc}}\in
C^{2+\alpha }\left( \mathbb{R}^{3}\setminus \Omega \right) .$ Hence, the
knowledge of the function $u_{\mathrm{sc}}(x,k)$ outside of the domain $%
\Omega $ yields the additional boundary data $g_{1}\left( x,k\right) ,$
where 
\begin{equation}
g_{1}(x,k)=\partial _{n}u(x,k),\quad x\in \partial \Omega ,k\in \lbrack 
\underline{k},\overline{k}].  \label{2.8}
\end{equation}

Even though we assume that the boundary measurements in (\ref{2.6}) are
conducted on the entire boundary $\partial \Omega ,$ this is done for the
analytical purpose only. In our computations we assume that we have only
backscattering data, which better suits our above mentioned target
application to imaging and identification of mine-like targets. We
complement the backscattering data on the rest of the boundary $\partial
\Omega $ by suitable values, see Section \ref{sec:8}.

Since we use here only a single direction of the propagation of the incident
plane wave $u_{0}\left( x,k\right) ,$ this is a problem with single
measurement data. All currently known uniqueness theorems for $n-$%
dimensional CIPs, $n\geq 2,$ with single measurement data are proven using
the method of \cite{BukhKlib}. This method is based on Carleman estimates.
Many publications of different authors have discussed this method. Since the
current work is not a survey of the technique of \cite{BukhKlib}, we refer
here only to a few such publications \cite{BK1,Im1,Klib92,Ksurvey,Trig,Yam}%
.\ In particular, \cite{Ksurvey} and \cite{Yam} are surveys of that method.
However, the technique of \cite{BukhKlib} works only if zero in the right
hand side of (\ref{2.3}) is replaced by such a function $f\in C\left( 
\overline{\Omega }\right) ,$ which does vanish in $\overline{\Omega }.$
Thus, since we study a numerical method here rather than the question of
uniqueness, we assume uniqueness of our CIP.

We model the propagation of the electric wave field in $\mathbb{R}^{3}$ by
the solution of the problem (\ref{2.3})-(\ref{2.5}). This modeling was
numerically justified in \cite{BMM}. It was demonstrated numerically in \cite%
{BMM} that this modeling can replace the modeling via the full Maxwell's
system, provided that only a single component of the electric field is
incident upon the medium. Then this component dominates two others and its
propagation is well governed by the time domain analog of equation (\ref{2.3}%
). This conclusion was verified via accurate imaging using electromagnetic
experimental data in, e.g. Chapter 5 of \cite{BK1} and \cite%
{KFB,KSNF1,TBKF1,TBKF2}.

\section{The asymptotic behavior of the function $u\left( x,k\right) $ as $k$
tends to infinity}

\label{sec:3}

To establish this asymptotic behavior, we use geodesic lines generated by
the function $c(x).$ Hence, we consider these lines in this section. The
discussion of this section is a modification of the corresponding discussion
of {\cite{KR}}. The Riemannian metric generated by the function $c(x)$ is 
\begin{equation*}
d\tau =\sqrt{c\left( x\right) }\left\vert dx\right\vert ,|dx|=\sqrt{%
(dx_{1})^{2}+(dx_{2})^{2}+(dx_{3})^{2}}.
\end{equation*}%
Consider the plane $P=\left\{ x_{3}=-R\right\} .$ Then $P\cap \overline{%
\Omega }=\varnothing .$ Consider unit vectors $e_{1}=\left( 1,0,0\right) $, $%
e_{2}=\left( 0,1,0\right) ,e_{3}=\left( 0,0,1\right) .$ An arbitrary point $%
\xi _{0}\in P$ can be represented as 
\begin{equation}
\xi _{0}=\xi _{0}(\eta _{1},\eta _{2})=\eta _{1}e_{1}+\eta
_{2}e_{2}-Re_{3},\quad (\eta _{1},\eta _{2})\in \mathbb{R}^{2}.  \label{3.1}
\end{equation}%
Let the function $\tau (x)$ be the solution of the following Cauchy problem
for the eikonal equation: 
\begin{equation}
\left\{ 
\begin{array}{c}
\left( \nabla \tau (x)\right) ^{2}=c(x), \\ 
\tau (x)=x_{3}\text{ for }x_{3}\leq -R.%
\end{array}%
\right.  \label{3.2}
\end{equation}
It is well known that $|\tau (x)|$ is the Riemannian distance between the
point $x$ and the plane $P$. Physically, $|\tau (x)|$ is the travel time
between the point $x$ and the plane $P$. To find the function $\tau (x)$ for 
$x_{3}>-R,$ it is necessary to solve the problem (\ref{3.2}). It is well
known that to solve this problem, one needs to solve a system of ordinary
differential equations. These equations also define geodesic lines of the
Riemannian metric. They are (see, e.g. \cite{R3}): 
\begin{equation}
\frac{d\xi }{ds}=\frac{p}{c(\xi )},\quad \frac{dp}{ds}=\frac{1}{2}\nabla
\left( \ln c(\xi )\right) ,\quad \frac{d\tau }{ds}=1,  \label{3.3}
\end{equation}%
where $s$ is a parameter and $p=\nabla \tau (\xi )$. Consider an arbitrary
point $\xi _{0}(\eta _{1},\eta _{2})\in P$ and the solution of the equations
(\ref{3.3}) with the Cauchy data 
\begin{equation}
\xi |_{s=0}=\xi _{0}(\eta _{1},\eta _{2}),\quad p|_{s=0}=\sqrt{c(\xi
_{0}(\eta _{1},\eta _{2}))}e_{3},\quad \tau |_{s=0}=0,  \label{3.4}
\end{equation}%
The solution of the problem (\ref{3.3}), (\ref{3.4}) defines the geodesic
line which passes through the point $\xi _{0}(\eta _{1},\eta _{2})$ in the
direction $e_{3}$. Hence, this line intersects the plane $P$ orthogonally.
For $s>0$ this solution determines the geodesic line $\xi =r_{1}(s,\eta
_{1},\eta _{2})$ and the vector $p=r_{2}(s,\eta _{1},\eta _{2}).$ This
vector lays in the tangent direction to that geodesic line. It is well known
from the theory of Ordinary Differential Equations that if the function $%
c(x)\in C^{m}(\mathbb{R}^{3})$, $m\geq 2$, then $r_{1}$ and $r_{2}$ are $%
C^{m-1}-$smooth functions.

By (\ref{3.1}) and (\ref{3.4}) 
\begin{equation}
\frac{\partial \xi }{\partial \eta _{1}}\Big|_{s=0}=e_{1},\quad \frac{%
\partial \xi }{\partial \eta _{2}}\Big|_{s=0}=e_{2}.  \label{3.5}
\end{equation}%
Noting that $c(\xi _{0}(\eta _{1},\eta _{2}))=1,$ we obtain from (\ref{3.4}) 
\begin{equation*}
\frac{d\xi }{ds}\Big|_{s=0}=\frac{\sqrt{c(\xi _{0}(\eta _{1},\eta _{2}))}}{%
c(\xi _{0}(\eta _{1},\eta _{2}))}e_{3}=e_{3}.
\end{equation*}%
Hence, 
\begin{equation}
\left\vert \frac{\partial (\xi _{1},\xi _{2},\xi _{3})}{\partial (s,\eta
_{1},\eta _{2})}\right\vert _{s=0}=1\neq 0.  \label{3.6}
\end{equation}

By (\ref{3.6}) the equality $x=r_{1}(s,\eta _{1},\eta _{2})$ can be uniquely
solved with respect to $s,\eta _{1},\eta _{2}$ for those points $x$ which
are sufficiently close to the plane $P$, as $s=s(x),\eta _{1}=\eta
_{1}(x),\eta _{2}=\eta _{2}(x).$ Hence, the equation 
\begin{equation*}
\xi =r_{1}(s,\eta _{1}(x),\eta _{2}(x))=\widehat{r}_{1}(s,x),\>s\in \left[
0,s(x)\right]
\end{equation*}%
defines the geodesic line $\Gamma (x)$ that passes through points $x$ and $%
\xi _{0}(\eta _{1}(x),\eta _{2}(x)):=\xi _{0}(x)$ and intersects the plane $%
P $ orthogonally.\ Extend the curve $\Gamma (x)$ for $x_{3}<-R$ as the
straight line by the equation $\xi =\xi _{0}(x)+se_{3}$, $s<0$. The
Riemannian distance between points $x$ and $\xi _{0}(x)$ is $s(x)=\tau (x)$.
Note that $\widehat{r}_{j}(s,x)=r_{j}(s,\eta _{1}(x),\eta _{2}(x)),j=1,2$
are $C^{m-1}-$smooth functions of their arguments. Since $\widehat{r}%
_{2}(s,x)=\nabla _{x}\tau (x)$, then $\tau (x)$ is the $C^{m}-$smooth
function. In our case $\tau (x)$ is $C^{15}-$smooth function of $x$ for
those points $x$ which are sufficiently close to the plane $P$.

We have constructed above the family of geodesic lines $\Gamma (x)$\ only
\textquotedblleft locally", i.e. only for those points $x$\ which are
located sufficiently close to the plane $P.$\ However, we need to consider
these lines \textquotedblleft globally". Hence, everywhere below we rely on
the following Assumption:

\textbf{Assumption 3.1}. \emph{We assume that above constructed geodesic
lines satisfy the regularity condition in }$\mathbb{R}^{3}$\emph{. In other
words, for each point }$x\in \mathbb{R}^{3}$ \emph{there exists a single
geodesic line }$\Gamma \left( x\right) $ \emph{connecting} $x$\emph{\ with
the plane }$P$ \emph{such that } $\Gamma \left( x\right) $ \emph{\
intersects }$P$\emph{\ orthogonally and the function $\tau(x)\in C^{15}(%
\mathbb{R}^3)$. }

A sufficient condition for the regularity of geodesic lines was derived in 
\cite{R4}, 
\begin{equation*}
\sum_{i,j=1}^{3}\frac{\partial ^{2}\ln c(x)}{\partial x_{i}\partial x_{j}}%
\xi _{i}\xi _{j}\geq 0,\>\forall x,\xi \in \mathbb{R}^{3}.
\end{equation*}%
Define the function $A(x)$ as 
\begin{equation}
A(x)=\left\{ 
\begin{array}{rl}
\exp \left( -\frac{1}{2}\int\limits_{\Gamma (x)}c^{-1}(\xi )\Delta _{\xi
}\tau (\xi )d\sigma \right) , & x_{3}>-R, \\ 
1, & x_{3}\leq -R.%
\end{array}%
\right.  \label{3.7}
\end{equation}%
The reason of the second line of (\ref{3.7}) is the second line of (\ref{3.2}%
) as well as the fact that $\Gamma (x)$ is the straight line for $x_{3}\leq
-R.$ Lemma 3.1 was proved in \cite{KR}, see Theorem 1 and the formula (4.25)
in this reference. In the proof of Theorem 1 of \cite{KR}, the $C^{15}-$%
smoothness of the function $c(x)$ was essentially used.

\textbf{Lemma 3.1. }\emph{Assume that conditions (\ref{2.0}) and (\ref{2.2})
are satisfied. Also, let Assumption 3.1 be in place. Then the following
asymptotic behavior of the solution }$u\left( x,k\right) $\emph{\ of the
problem (\ref{2.3})-(\ref{2.5}) holds: }%
\begin{equation}
u\left( x,k\right) =A(x)e^{-ik\tau \left( x\right) }\left( 1+O\left(
1/k\right) \right) ,k\rightarrow \infty ,x\in \Omega .  \label{3.8}
\end{equation}%
\emph{Here }$\left\vert O\left( 1/k\right) \right\vert \leq B_{1}/k,\forall
x\in \overline{\Omega },$\emph{\ where the constant }$B_{1}=B_{1}\left(
\Omega ,c\right) >0$\emph{\ depends only on listed parameters.}

Hence, it follows from this lemma and (\ref{3.7}) that there exists a number 
$k_{0}\left( \Omega ,c\right) >0$ depending only on listed parameters such
that%
\begin{equation}
\left\vert O\left( 1/k\right) \right\vert <1/2,\forall k\geq k_{0}\left(
\Omega ,c\right) ,\forall x\in \overline{\Omega }.  \label{3.80}
\end{equation}%
\begin{equation}
u\left( x,k\right) \neq 0,\forall x\in \overline{\Omega },\forall k\geq
k_{0}\left( \Omega ,c\right) .  \label{3.9}
\end{equation}

\section{Using the Lippmann-Schwinger equation}

\label{sec:4}

In this section we use the Lippmann-Schwinger equation to derive some
important facts, which we need both for our algorithm and for the
convergence analysis. We are essentially using here results of Chapter 8 of
the book of Colton and Kress \cite{CK}. In accordance with the
regularization theory, we need to assume that there exists unique exact
solution $c^{\ast }$ of our CIP for the noiseless data $g^{\ast }\left(
x,k\right) $ in (\ref{2.6}) \cite{BK1,T}. Everywhere below the superscript
\textquotedblleft $\ast $" denotes functions generated by $c^{\ast }.$

Denote 
\begin{equation*}
\Phi \left( x,y\right) =\frac{\exp \left( -i\overline{k}\left\vert
x-y\right\vert \right) }{4\pi \left\vert x-y\right\vert },x\neq y.
\end{equation*}%
In this section the function $\beta \in C^{\alpha }\left( \overline{\Omega }%
\right) $ and satisfies condition (\ref{2.2}). The Lippmann-Schwinger
equation for the function $u\left( x\right) :=u\left( x,\overline{k}\right) $
is%
\begin{equation}
u\left( x\right) =\exp \left( -i\overline{k}x_{3}\right) +\overline{k}%
^{2}\dint\limits_{\Omega }\Phi \left( x,y\right) \beta (y)u\left( y\right)
dy.  \label{8.3}
\end{equation}%
If the function $u\left( x\right) $ satisfies equation (\ref{8.3}) for $x\in
\Omega ,$ then we can extend it for $x\in \mathbb{R}^{3}\diagdown \Omega $
via substitution these points $x$ in the right hand side of (\ref{8.3}).
Hence, to solve (\ref{8.3}), it is sufficient to find the function $u\left(
x\right) $ only for points $x\in \Omega .$ Consider the linear operator $%
K_{\beta }$ defined as 
\begin{equation}
\left( K_{\beta }u\right) \left( x\right) =\overline{k}^{2}\dint\limits_{{%
\Omega }}\Phi \left( x,y\right) \beta (y)u\left( y\right) dy.  \label{8.30}
\end{equation}%
It follows from Theorem 8.1 of \cite{CK} that 
\begin{equation}
K_{\beta }:C^{\alpha }\left( \overline{\Omega }\right) \rightarrow
C^{2+\alpha }\left( \mathbb{R}^{3}\right) \text{ and }\left\Vert K_{\beta
}u\right\Vert _{2+\alpha }\leq B_{2}|\beta |_{\alpha }\left\vert
u\right\vert _{\alpha }.  \label{8.4}
\end{equation}%
Here and below $B_{2}=B_{2}\left( \beta ,\overline{k},\Omega _{1},\Omega
\right) >0$ denotes different constants which depend only on listed
parameters. Therefore, the operator $K_{\beta }$ maps $C^{\alpha }\left( 
\overline{\Omega }\right) $ in $C^{\alpha }\left( \overline{\Omega }\right) $
as a compact operator, $K_{\beta }:C^{\alpha }\left( \overline{\Omega }%
\right) \rightarrow C^{\alpha }\left( \overline{\Omega }\right) $. Hence,
the Fredholm theory is applicable to equation (\ref{8.3}). Lemmata 4.1 and
4.2 follow from Theorem 8.3 and Theorem 8.7 of \cite{CK} respectively.

\textbf{Lemma 4.1}. \emph{The function }$u\in C^{2+\alpha }\left( \mathbb{R}%
^{3}\right) $\emph{\ is a solution of the problem (\ref{2.3})-(\ref{2.5}) if
and only if it is a solution of equation (\ref{8.3}).}

\textbf{Lemma 4.2}. \emph{There exists unique solution }$u\in C^{2+\alpha
}\left( \mathbb{R}^{3}\right) $\emph{\ of the problem (\ref{2.3})-(\ref{2.5}%
). Consequently (Lemma 4.1) there exists unique solution }$u\in C^{2+\alpha
}\left( \mathbb{R}^{3}\right) $\emph{\ of the problem (\ref{8.3}) and these
two solutions coincide. Furthermore, by the Fredholm theory }$\left\vert
u\right\vert _{\alpha }\leq B_{2}\left\vert \exp \left( -i\overline{k}%
x_{3}\right) \right\vert _{\alpha }.$ \emph{Also, with a different constant }%
$B_{2}$,\emph{\ }$\left\vert u\right\vert _{2+\alpha }\leq B_{2}\left\vert
\exp \left( -i\overline{k}x_{3}\right) \right\vert _{2+\alpha }.$\emph{\ }

\textbf{Proof}. We need to prove only the last estimate. By (\ref{8.4}) 
\begin{equation}
\left\vert K_{\beta }u\right\vert _{2+\alpha }\leq B_{2}|\beta |_{\alpha
}\left\vert u\right\vert _{\alpha }\leq B_{2}\left\vert \exp \left( -i%
\overline{k}x_{3}\right) \right\vert _{\alpha }\leq B_{2}\left\vert \exp
\left( -i\overline{k}x_{3}\right) \right\vert _{2+\alpha }.  \label{8.40}
\end{equation}%
On the other hand, by (\ref{8.3}) 
\begin{equation}
\left\vert u\right\vert _{2+\alpha }\leq \left\vert \exp \left( -i\overline{k%
}x_{3}\right) \right\vert _{2+\alpha }+\left\vert K_{\beta }u\right\vert
_{2+\alpha }.  \label{8.41}
\end{equation}%
Thus, the desired estimate follows from (\ref{8.40}) and (\ref{8.41}). \ $%
\square $

Lemma 4.3 follows from Lemmata 4.1, 4.2 and results of Chapter 9 of the book
of Vainberg \cite{V}.

\textbf{Lemma 4.3}. \emph{For all }$x\in \overline{\Omega },k>0$\emph{\ the
function }$u\left( x,k\right) $ \emph{is infinitely many times
differentiable with respect to }$k$\emph{.\ Furthermore, }$\partial
_{k}^{n}u\in C^{2+\alpha }\left( \overline{\Omega }\right) $\emph{\ and }%
\begin{equation*}
\lim_{\omega \rightarrow 0,\omega \in \mathbb{R}}\left\vert \partial
_{k}^{n}u\left( x,k+\omega \right) -\partial _{k}^{n}u\left( x,k\right)
\right\vert _{2+\alpha }=0,\emph{\ }n=0,1,...\emph{\ }
\end{equation*}

Let the function $\chi \in C^{2}\left( \mathbb{R}^{3}\right) $ be such that%
\begin{equation}
\chi \left( x\right) =\left\{ 
\begin{array}{c}
1,\text{ if }x\in \Omega _{1}, \\ 
\text{between }0\text{ and }1,\text{ if }x\in \Omega \diagdown \Omega _{1},
\\ 
0,\text{ if }x\in \mathbb{R}^{3}\diagdown \Omega _{1}.%
\end{array}%
\right.  \label{8.60}
\end{equation}
The existence of such functions is well known from the Real Analysis course.
Consider a complex valued function $\rho \left( x\right) \in C^{\alpha
}\left( \overline{\Omega }\right) .$ Let $\widehat{\rho }\left( x\right)
=\chi \left( x\right) \rho \left( x\right) .$ Then 
\begin{equation}
\widehat{\rho }\in C^{\alpha }\left( \mathbb{R}^{3}\right) ,\widehat{\rho }%
\left( x\right) =0\text{ in }\mathbb{R}^{3}\diagdown \Omega .  \label{8.7}
\end{equation}

\textbf{Theorem 4.1}. \emph{Assume that the exact coefficient }$c^{\ast
}\left( x\right) ,$\emph{\ satisfies conditions (\ref{2.0}), (\ref{2.2}).
Let }$\beta ^{\ast }\left( x\right) =c^{\ast }\left( x\right) -1$\emph{. Let 
}$u^{\ast }\left( x,\overline{k}\right) $\emph{\ be the solution of the
problem (\ref{2.3})-(\ref{2.5}) in which }$\beta \left( x\right) $\emph{\ is
replaced with }$\beta ^{\ast }\left( x\right) .$\emph{\ Consider equation (%
\ref{8.3}), in which }$\beta \left( x\right) $\emph{\ is replaced with }$%
\widehat{\rho }\left( x\right) ,$%
\begin{equation}
u_{\rho }\left( x,\overline{k}\right) =\exp \left( -i\overline{k}%
x_{3}\right) +\overline{k}^{2}\dint\limits_{\Omega }\Phi \left( x,y\right) 
\widehat{\rho }(y)u_{\rho }\left( y,\overline{k}\right) dy,x\in \Omega .
\label{8.8}
\end{equation}%
\emph{\ Then there exists a sufficiently small number }$\theta ^{\ast
}=\theta ^{\ast }\left( \beta ^{\ast },\overline{k},\chi ,\Omega _{1},\Omega
\right) \in \left( 0,1\right) $\emph{\ depending only on listed parameters
such that if }$\left\vert \rho -\beta ^{\ast }\right\vert _{\alpha }\leq
\theta $ \emph{and }$\theta \in \left( 0,\theta ^{\ast }\right) ,$\emph{\
then equation (\ref{8.8}) has unique solution }$u_{\rho }\in C^{\alpha
}\left( \overline{\Omega }\right) .$ \emph{Furthermore, the function }$%
u_{\rho }\in C^{2+\alpha }\left( \mathbb{R}^{3}\right) $\emph{\ and }%
\begin{equation}
\left\vert u_{\rho }\left( x,\overline{k}\right) -u^{\ast }\left( x,%
\overline{k}\right) \right\vert _{2+\alpha }\leq Z^{\ast }\theta ,
\label{8.9}
\end{equation}%
\emph{\ where the constant }$Z^{\ast }=Z^{\ast }\left( \beta ^{\ast },%
\overline{k},\chi ,\Omega _{1},\Omega \right) >0$\emph{\ depends only on
listed parameters.}

\textbf{Proof}. Below $Z^{\ast }$ denotes different positive constants
depending on the above parameters. We have $\beta ^{\ast }\left( x\right)
=\chi \left( x\right) \beta ^{\ast }\left( x\right) +\left( 1-\chi \left(
x\right) \right) \beta ^{\ast }\left( x\right) .$ Since by (\ref{2.2}) the
function $\beta ^{\ast }\left( x\right) =0$ outside of the domain $\Omega
_{1}$, then (\ref{8.60}) implies that $\left( 1-\chi \left( x\right) \right)
\beta ^{\ast }\left( x\right) \equiv 0.$ Hence, $\beta ^{\ast }\left(
x\right) =\chi \left( x\right) \beta ^{\ast }\left( x\right) .$ Hence, $%
\left( \widehat{\rho }-\beta ^{\ast }\right) \left( x\right) =\chi \left(
x\right) \left( \rho -\beta ^{\ast }\right) \left( x\right) .$ Hence, using
notation (\ref{8.30}), we rewrite equation (\ref{8.8}) in the following
equivalent form:%
\begin{equation}
\left( I-K_{\beta ^{\ast }}\right) u_{\rho }=\exp \left( -i\overline{k}%
x_{3}\right) +A_{\rho -\beta ^{\ast }}\left( u_{\rho }\right) ,x\in \Omega ,
\label{8.10}
\end{equation}%
\begin{equation}
A_{\rho -\beta ^{\ast }}\left( u_{\rho }\right) \left( x\right) =\overline{k}%
^{2}\dint\limits_{\Omega }\Phi \left( x,y\right) \chi \left( y\right) \left(
\rho -\beta ^{\ast }\right) (y)u_{\rho }\left( y,\overline{k}\right) dy,x\in
\Omega .  \label{8.100}
\end{equation}%
The linear operator $A_{\rho -\beta ^{\ast }}:C^{\alpha }\left( \overline{%
\Omega }\right) \rightarrow C^{\alpha }\left( \overline{\Omega }\right) .$
We have for any function $p\in C^{\alpha }\left( \overline{\Omega }\right) $ 
\begin{equation}
\left\vert A_{\rho -\beta ^{\ast }}\left( p\right) \right\vert _{\alpha
}=\left\vert \overline{k}^{2}\dint\limits_{\Omega }\Phi \left( x,y\right)
\chi \left( y\right) \left( \rho -\beta ^{\ast }\right) (y)p\left( y\right)
dy\right\vert _{2+\alpha }\leq Z^{\ast }\theta |p|_{\alpha }.  \label{8.101}
\end{equation}%
Therefore, 
\begin{equation}
\left\Vert A_{\rho -\beta ^{\ast }}\right\Vert \leq Z^{\ast }\theta .
\label{8.11}
\end{equation}%
It follows from Lemma 4.2 and the Fredholm theory that the operator $\left(
I-K_{\beta ^{\ast }}\right) $ has a bounded inverse operator $T=\left(
I-K_{\beta ^{\ast }}\right) ^{-1},T:C^{\alpha }\left( \overline{\Omega }%
\right) \rightarrow C^{\alpha }\left( \overline{\Omega }\right) $ and also 
\begin{equation}
\left\Vert T\right\Vert \leq Z^{\ast }.  \label{8.111}
\end{equation}
Hence, using (\ref{8.10}), we obtain%
\begin{equation}
u_{\rho }=T\left( \exp \left( -i\overline{k}x_{3}\right) \right) +\left(
TA_{\rho -\beta ^{\ast }}\right) \left( u_{\rho }\right) .  \label{8.12}
\end{equation}%
It follows from (\ref{8.11}) and (\ref{8.111}) that there exists a
sufficiently small number $\theta ^{\ast }\in (0,1)$ depending on $\beta
^{\ast }$, $\overline{k},$ $\chi ,\Omega _{1},\Omega $ such that if $\theta
\in \left( 0,\theta ^{\ast }\right) $ and $|\rho -\beta ^{\ast }|_{\alpha
}\leq \theta ,$ then the operator $\left( TA_{\rho -\beta ^{\ast }}\right)
:C^{\alpha }\left( \overline{\Omega }\right) \rightarrow C^{\alpha }\left( 
\overline{\Omega }\right) $ is contraction mapping. This implies uniqueness
and existence of the solution $u_{\rho }\in C^{\alpha }\left( \overline{%
\Omega }\right) $ of equation (\ref{8.12}), which is equivalent to equation (%
\ref{8.8}). Also, 
\begin{equation}
\left\vert u_{\rho }\right\vert _{\alpha }\leq Z^{\ast }\left\vert \exp
\left( -i\overline{k}x_{3}\right) \right\vert _{\alpha },  \label{8.13}
\end{equation}%
with a different constant $Z^{\ast }.$ Furthermore, by (\ref{8.4}) the
function $u_{\rho }\in C^{2+\alpha }\left( \overline{\Omega }\right) .$

We now prove estimate (\ref{8.9}). Let $\widetilde{u}\left( x\right)
=u^{\ast }\left( x,\overline{k}\right) -u_{\rho }\left( x,\overline{k}%
\right) .$ Since $\left( I-K_{\beta ^{\ast }}\right) u^{\ast }=\exp \left( -i%
\overline{k}x_{3}\right) ,$ we obtain the following analog of (\ref{8.10}) 
\begin{equation}
\left( I-K_{\beta ^{\ast }}\right) \widetilde{u}=A_{\rho -\beta ^{\ast
}}\left( u_{_{\rho }}\right) ,x\in \Omega .  \label{8.14}
\end{equation}%
Hence, $\widetilde{u}=\left( TA_{\rho -\beta ^{\ast }}\right) \left( u_{\rho
}\right) .$ Hence, (\ref{8.11}) and (\ref{8.13}) lead to 
\begin{equation}
\left\vert \widetilde{u}\right\vert _{\alpha }\leq Z^{\ast }\theta .
\label{8.15}
\end{equation}%
Next, we rewrite (\ref{8.14}) as%
\begin{equation}
\widetilde{u}=K_{\beta ^{\ast }}\widetilde{u}+A_{\rho -\beta ^{\ast }}\left(
u_{\rho }\right) ,\text{ }x\in \Omega .  \label{8.16}
\end{equation}%
By (\ref{8.4}) and (\ref{8.100}) the right hand side of equation (\ref{8.16}%
) belongs to the space $C^{2+\alpha }\left( \overline{\Omega }\right) .$
Hence, using (\ref{8.4}), (\ref{8.101}), (\ref{8.13}) and (\ref{8.15}), we
obtain from (\ref{8.16}) that $\left\vert \widetilde{u}\right\vert
_{2+\alpha }\leq Z^{\ast }\theta .$ $\square $

\section{Numerical method}

\label{sec:5}

\subsection{Some auxiliary functions}

\label{subsection 5.1}

Starting from this subsection and until section 8 $x\in \Omega ,$ and we do
not consider $x\in \mathbb{R}^{3}\diagdown \Omega .$ We define in this
subsection the logarithm of the complex valued function $u(x,k)$, $x\in
\Omega ,k>0$. We note that, except of subsection 5.3, we use only
derivatives of $\log u$ and do not use $\log u$ itself$.$ Since $\partial
\log u=\partial u/u,$ then this eliminates the uncertainty linked with $%
\func{Im}\left( \log u\right) .$ Below $\overline{k}>\underline{k}%
>k_{0}\left( \Omega ,c\right) $ and we consider $k\in \left[ \underline{k},%
\overline{k}\right] .$ The number $k_{0}\left( \Omega ,c\right) $ was
defined in (\ref{3.80}), (\ref{3.9}). Hence, by (\ref{3.9}) $u\left(
x,k\right) \neq 0$ for $k\in \left[ \underline{k},\overline{k}\right] ,x\in 
\overline{\Omega }.$ It is convenient to consider in this section only $%
C^{2+\alpha }-$smoothness of the function $u\left( x,k\right) ,$ i.e. $u\in
C^{2+\alpha }\left( \overline{\Omega }\right) $.

By a simple calculation, $\mathrm{curl}\left( \nabla u(x,\overline{k})/u(x,%
\overline{k})\right) =0$ in $\mathbb{R}^{3}$. Since $\Omega $ is a convex
domain, then there exists a function $V\left( x\right) $, such that 
\begin{equation}
\nabla V\left( x\right) =\frac{\nabla u\left( x,\overline{k}\right) }{%
u\left( x,\overline{k}\right) },x\in \overline{\Omega }.  \label{4.0}
\end{equation}%
By (\ref{4.0}) 
\begin{equation*}
e^{-V(x)}\left( u(x,\overline{k})\nabla V(x)-\nabla u(x,\overline{k})\right)
=0,x\in \overline{\Omega },
\end{equation*}%
which implies $\nabla (e^{-V(x)}u(x,\overline{k}))=0$. Thus, there exists a
constant $C$ such that $u(x,\overline{k})=Ce^{V(x)},x\in \overline{\Omega }.$
Since the function $V(x)$ is uniquely determined up to an addition of a
constant, we can choose $V(x)$ such that $C=1$. In summary, we can find a
function $V$ such that 
\begin{equation}
u(x,\overline{k})=e^{V(x)},\mbox{or }V(x)=\log (u(x,\overline{k})).
\label{4.5}
\end{equation}%
Since the function $u\in C^{2+\alpha }\left( \overline{\Omega }\right) ,$
then it follows from (\ref{4.0}) that $\nabla V\in C^{1+\alpha }\left( 
\overline{\Omega }\right) $.

By Lemma 4.3 the derivative $\partial _{k}u\left( x,k\right) \in C^{2+\alpha
}\left( \overline{\Omega }\right) $ exists and the continuity property (\ref%
{10}) is valid. Hence, we can define the function $v(x,k)$ for all $k\in %
\left[ \underline{k},\overline{k}\right] ,x\in \overline{\Omega }$ as 
\begin{equation}
v(x,k)=-\int_{k}^{\overline{k}}\frac{\partial _{k}u(x,\kappa )}{u(x,\kappa )}%
d\kappa +V(x),x\in \overline{\Omega }.  \label{def vxk}
\end{equation}%
Differentiate (\ref{def vxk}) with respect to $k$. We obtain $\partial
_{k}v(x,k)=\partial _{k}u(x,k)/u(x,k)$. Therefore, 
\begin{equation*}
e^{-v(x,k)}\left( u(x,k)\partial _{k}v(x,k)-\partial _{k}u(x,k)\right)
=0,x\in \overline{\Omega },
\end{equation*}%
which implies $\partial _{k}(e^{-v(x,k)}u(x,k))=0$ or $u(x,k)=C(x)e^{v(x,k)}$
for all $k\in \left[ k_{0}\left( \Omega ,c\right) ,\overline{k}\right] ,x\in 
\overline{\Omega }.$ In particular, taking $k=\overline{k}$ and using %
\eqref{4.5}, we obtain $C\left( x\right) =1.$

\begin{lemma}
\label{Lemma 5.1111} \emph{For each }$k\in \left[ \underline{k},\overline{k}%
\right] $\emph{\ the gradient } $\nabla v\in C^{1+\alpha }\left( \overline{%
\Omega }\right) .$\emph{\ In addition for all } $x\in \overline{\Omega }:$

\begin{enumerate}
\item \label{item 1} $u(x, k) = e^{v(x, k)}$,

\item \label{item 2} $\nabla v(x, k) = \nabla u(x, k)/u(x, k),$

\item \label{item 3} $\partial_k v(x, k) = \displaystyle \frac{\partial_k
u(x, k)}{u(x, k)}$,

\item \label{item 4} 
\begin{equation}  \label{4.1}
\Delta v(x, k) + (\nabla v(x, k))^2 = -k^2 c(x).
\end{equation}
\end{enumerate}
\end{lemma}

\textbf{Proof}. The smoothness of $\nabla v$ follows from $\nabla V\in
C^{1+\alpha }\left( \overline{\Omega }\right) $ and from \eqref{def vxk}.
Item 1 was established above. The differentiation of the equality of item 1
leads to item 2. Item 3 follows from \eqref{def vxk}. Equation \eqref{4.1}
follows from item 1 and \eqref{2.3}. $\square$

Thus, $V\left( x\right) $ is our tail function. The exact tail function,
which is generated by the exact coefficient $c^{\ast }\left( x\right) ,$ is $%
V^{\ast }\left( x\right) =\log u^{\ast }\left( x,\overline{k}\right) .$

\subsection{Integral differential equation}

\label{subsection 5.2}

Consider the function $q\left( x,k\right) $ defined as 
\begin{equation}
q(x,k)=\partial _{k}v(x,k)=\frac{\partial _{k}u(x,k)}{u(x,k)},\quad x\in
\Omega ,k\in \left( \underline{k},\overline{k}\right) .  \label{def q}
\end{equation}%
Here, we have used item \ref{item 3} in Lemma \ref{Lemma 5.1111} for the
latter fact. By \eqref{def vxk} 
\begin{equation}
v(x,k)=-\int_{k}^{\overline{k}}q(x,\kappa )d\kappa +V(x),\quad x\in \Omega
,k\in (\underline{k},\overline{k}).  \label{4.4}
\end{equation}%
Note that $V\left( x\right) =v\left( x,\overline{k}\right) .$ We call $%
V\left( x\right) $ the \textquotedblleft tail function". We note that the
number $\overline{k}$ plays the role of the regularization parameter of our
numerical method. A convergence analysis of our method for $\overline{k}%
\rightarrow \infty $ is a very challenging problem and we do not yet know
how to address it. The differentiation of \eqref{4.1} with respect to $k$
leads to 
\begin{equation*}
\Delta q(x,k)+2\nabla q(x,k)\nabla v(x,k)=-2kc(x)=2(\Delta v+(\nabla
v)^{2})/k.
\end{equation*}%
This, \ and (\ref{4.4}) imply that for all $k\in \left[ \underline{k},%
\overline{k}\right] $%
\begin{eqnarray}
&&k\Delta q(x,k)+2k\nabla q(x,k)\nabla \left( -\int_{k}^{\overline{k}%
}q(x,\kappa )d\kappa +V(x)\right)  \notag \\
&=&2\left( \Delta \left( -\int_{k}^{\overline{k}}q(x,\kappa )d\kappa
+V(x)\right) +\left( \nabla \left( -\int_{k}^{\overline{k}}q(x,\kappa
)d\kappa +V(x)\right) \right) ^{2}\right) .  \label{4.6}
\end{eqnarray}%
By Lemma \ref{Lemma 5.1111} as well as by (\ref{2.6}) and \eqref{def q}, the
function $q$ satisfies the Dirichlet boundary condition 
\begin{equation}
q(x,k)=\frac{\partial _{k}g(x,k)}{g(x,k)}=:\psi (x,k)\quad \mbox{on }%
\partial \Omega ,\psi \left( x,k\right) \in C^{2+\alpha }\left( \partial
\Omega \right) ,\forall k\in \left[ \underline{k},\overline{k}\right] .
\label{4.7}
\end{equation}

We have obtained a nonlinear integral differential equation (\ref{4.6}) for
the function $q\left( x,k\right) $ with the Dirichlet boundary condition (%
\ref{4.7}). Both functions $q$ and $V$ in (\ref{4.6}) are unknown.\ To solve
our inverse problem, both these functions need to be approximated. Here is a
brief description how do we do this. We start from finding a first
approximation $V_{0}\left( x\right) $ for the tail function, see subsection
5.3. To approximate the function $q$, we iteratively solve the problem (\ref%
{4.6}), (\ref{4.7}) inside of the domain $\Omega .$ Given an approximation
for $q$, we find the next approximation for the unknown coefficient $c,$ and
then we solve the Lippman-Schwinger equation inside of the domain $\Omega $
with this updated coefficient $c$. Next, we find the new approximation for
the gradient of the tail function $V$ via (\ref{4.0}) as $\nabla V=\nabla
u\left( x,\overline{k}\right) /u\left( x,\overline{k}\right) $ and similarly
for the new approximation for $\Delta V=\func{div}\left( \nabla V\right) $.
So, this is an analog of the well known predictor-corrector procedure, where
updates for $V$ are predictors and updates for $q$ and $c$ are correctors.

\subsection{The first approximation $V_{0}\left( x\right) $ for the tail
function}

\label{subsection 5.3}

Consider the exact coefficient $c^{\ast }\left( x\right) $ and assume that
conditions (\ref{2.0}), (\ref{2.2}) as well as Assumption 3.1 hold for $%
c^{\ast }\left( x\right) $. Then (\ref{3.8}) holds for $c\left( x\right)
:=c^{\ast }\left( x\right) .$ Assume that the number $\overline{k}$ is
sufficiently large. For all $k\geq \overline{k}$ drop the term $O\left(
1/k\right) $ in (\ref{3.8}). Hence, we approximate the function $u^{\ast
}\left( x,k\right) $ as%
\begin{equation}
u^{\ast }\left( x,k\right) =A^{\ast }\left( x\right) e^{-ik\tau ^{\ast
}\left( x\right) },k\geq \overline{k}.  \label{4.8}
\end{equation}%
We set 
\begin{equation*}
\log u^{\ast }\left( x,k\right) =\ln A^{\ast }\left( x\right) -ik\tau ^{\ast
}\left( x\right) \text{ for }k\geq \overline{k}.
\end{equation*}%
Hence, 
\begin{equation}
\log u^{\ast }\left( x,k\right) =-ik\tau ^{\ast }\left( x\right) \left(
1+O\left( \frac{1}{k}\right) \right) ,k\rightarrow \infty .  \label{101}
\end{equation}%
Drop again the term $O\left( 1/k\right) $ in (\ref{101}). Next, set $k=%
\overline{k}.$ Hence, we approximate the exact tail function $V^{\ast
}\left( x\right) $ for $k=\overline{k}$ as 
\begin{equation}
V^{\ast }\left( x\right) =-i\overline{k}\tau ^{\ast }\left( x\right) .
\label{102}
\end{equation}%
Using (\ref{def q}) and (\ref{4.8}), we obtain 
\begin{equation}
q^{\ast }\left( x,\overline{k}\right) =-i\tau ^{\ast }\left( x\right) .
\label{103}
\end{equation}%
Set in equation (\ref{4.6}) $k:=\overline{k},q\left( x,\overline{k}\right)
:=q^{\ast }\left( x,\overline{k}\right) ,V(x):=V^{\ast }(x).$ Next,
substitute in the resulting equation formulae (\ref{102}) and (\ref{103}).
Also, use (\ref{4.7}) for $\psi :=\psi ^{\ast }$. We obtain%
\begin{equation}
\begin{array}{c}
\Delta \tau ^{\ast }=0\text{ in }\Omega , \\ 
\tau ^{\ast }\mid _{\partial \Omega }=i\psi ^{\ast }\left( x,\overline{k}%
\right) .%
\end{array}
\label{4.110}
\end{equation}%
Thus, we have obtained the Dirichlet boundary value problem (\ref{4.110})
for the Laplace equation with respect to the function $\tau ^{\ast }\left(
x\right) $. Recalling that $\partial \Omega \in C^{2+\alpha }$ and that the
function $\psi ^{\ast }\left( x,\overline{k}\right) \in C^{2+\alpha }\left(
\partial \Omega \right) $ and applying the Schauder theorem \cite{Lad}, we
obtain that there exists unique solution $\tau ^{\ast }\in C^{2+\alpha
}\left( \overline{\Omega }\right) $ of the problem (\ref{4.110}).

In practice, however, we have the non-exact boundary data $\psi \left(
x,k\right) $ rather than the exact data $\psi ^{\ast }\left( x,k\right) .$
Thus, we set the first approximation $V_{0}\left( x\right) $ for the tail
function $V\left( x\right) $ as 
\begin{equation}
V_{0}\left( x\right) =-i\overline{k}\tau \left( x\right) ,  \label{106}
\end{equation}%
where the function $\tau \left( x\right) $ is the $C^{2+\alpha }\left( 
\overline{\Omega }\right) -$solution of the following analog of is the
solution of the problem (\ref{4.110}):%
\begin{equation}
\begin{array}{c}
\Delta \tau =0\text{ in }\Omega , \\ 
\tau \mid _{\partial \Omega }=i\psi \left( x,\overline{k}\right) .%
\end{array}
\label{1060}
\end{equation}

Theorem 5.1 estimates the difference between functions $V_{0}\left( x\right) 
$ and $V^{\ast }\left( x\right) .$

\textbf{Theorem 5.1}. \emph{Assume that relations (\ref{102}), (\ref{103}), (%
\ref{106}) and (\ref{1060}) are valid. Then there exists a constant }$%
C=C\left( \Omega \right) >0$\emph{\ depending only on the domain }$\Omega $%
\emph{\ such that}%
\begin{equation}
\left\vert V_{0}-V^{\ast }\right\vert _{2+\alpha }\leq C\overline{k}%
\left\Vert \psi \left( x,\overline{k}\right) -\psi ^{\ast }\left( x,%
\overline{k}\right) \right\Vert _{C^{2+\alpha }\left( \partial \Omega
\right) }.  \label{107}
\end{equation}

\textbf{Proof}. Note that (\ref{4.110}) follows from (\ref{4.6}), (\ref{4.7}%
), (\ref{102}) and (\ref{103}). Denote $\widetilde{\tau }\left( x\right)
=\tau \left( x\right) -\tau ^{\ast }\left( x\right) .$ Then (\ref{4.110})
and (\ref{1060}) imply that 
\begin{equation*}
\begin{array}{c}
\Delta \widetilde{\tau }=0,x\in \Omega , \\ 
\widetilde{\tau }\mid _{\partial \Omega }=i\left( \psi -\psi ^{\ast }\right)
\left( x,\overline{k}\right) .%
\end{array}%
\end{equation*}%
Hence, the Schauder theorem \cite{Lad} leads to (\ref{107}). $\square $

\textbf{Remarks 5.1}:

\begin{enumerate}
\item Theorem 5.1 means that the accuracy of the approximation of the exact
tail function $V^{\ast }$ by the function $V_{0}$ depends only on the
accuracy of the approximation of the exact boundary condition $\psi ^{\ast
}\left( x,\overline{k}\right) $ by the boundary condition $\psi \left( x,%
\overline{k}\right) .$ Let $\delta >0$ be the level of the error in the
boundary data at $k:=\overline{k}$, i.e. $\left\Vert \psi \left( x,\overline{%
k}\right) -\psi ^{\ast }\left( x,\overline{k}\right) \right\Vert
_{C^{2+\alpha }\left( \partial \Omega \right) }\leq \delta .$ Hence, if $%
\delta $ is sufficiently small, then, by (\ref{107}), the norm $\left\vert
V_{0}-V^{\ast }\right\vert _{2+\alpha }$ is also sufficiently small. In the
regularization theory the error in the data $\delta $ is always assumed to
be sufficiently small \cite{BK1,T}. Thus, we have obtained the tail function 
$V_{0}$ in a sufficiently small neighborhood of the exact tail function $%
V^{\ast }$. Furthermore, in doing so, we have not used any \emph{a priori}
knowledge about a sufficiently small neighborhood of the function $V^{\ast
}. $ The smallness of that neighborhood depends only on the level of the
error in the boundary data. The latter is exactly what is required in the
regularization theory.

\item Thus, the global convergence property is achieved just at the start of
our iterative process. It is achieved due to two factors. The first factor
is the elimination of the unknown coefficient from equation (\ref{4.1}) and
obtaining the integral differential equation (\ref{4.6}). The second factor
is dropping the term $O\left( 1/k\right) $ in\ (\ref{4.8}) and (\ref{101}).

\item Still, our numerical experience shows that we need to do more
iterations to obtain better accuracy. These iterations are described in
subsection 5.4.

\item We point out that we use the approximations (\ref{4.8}), (\ref{102})
and (\ref{103}) of the exact tail function $V^{\ast }\left( x\right) $ only
on the first iteration of our method: to obtain the first approximation $%
V_{0}\left( x\right) $ for the tail function. However, we do not use them on
follow up iterations. On a deeper sense, these approximations are introduced
because the problem of constructing globally convergent numerical methods
for CIPs is well known to be a \emph{tremendously challenging} one. Indeed,
CIPs are both nonlinear and ill-posed. Thus, it makes sense to use such an
approximation. Because of this approximation, one can also call our
technique an \emph{approximately globally convergent numerical method}, see
section 1.1.2 in \cite{BK1} as well as \cite{KSNF1} for detailed discussions
of the notion of the \emph{approximate global convergence}.
\end{enumerate}

\subsection{The algorithm}

\label{sec:5.4}

Let $h>0$ be the partition step size of a uniform partition of the frequency
interval $\left[ \underline{k},\overline{k}\right] $, 
\begin{equation}
\underline{k}=k_{N}<k_{N-1}<...<k_{1}<k_{0}=\overline{k},k_{j-1}-k_{j}=h.
\label{4.00}
\end{equation}%
Approximate the function $q\left( x,k\right) $ as a piecewise constant
function with respect to $k\in \left[ \underline{k},\overline{k}\right] .$
Then (\ref{4.7}) implies that the boundary condition $\psi \left( x,k\right)
,x\in \partial \Omega $ should also be approximated by a piecewise constant
function with respect to $k\in \left[ \underline{k},\overline{k}\right] .$
Let 
\begin{equation}
q\left( x,k\right) =q_{n}\left( x\right) ,\psi \left( x,k\right) =\psi
_{n}\left( x\right) ,k\in \left[ k_{n},k_{n-1}\right) ,n=1,...,N.
\label{4.16}
\end{equation}%
We set $q_{0}\left( x\right) \equiv 0.$ Denote%
\begin{equation}
\overline{q_{n-1}}=\dsum\limits_{j=0}^{n-1}q_{j}\left( x\right) .
\label{4.17}
\end{equation}%
Hence, (\ref{4.4}) becomes%
\begin{equation}
v\left( x,k\right) =-\left( k_{n-1}-k\right) q_{n}\left( x\right) -h%
\overline{q_{n-1}}+V\left( x\right) ,k\in \left[ k_{n},k_{n-1}\right) .
\label{4.18}
\end{equation}%
Hence, the problem (\ref{4.6}), (\ref{4.7}) can be rewritten for $k\in \left[
k_{n},k_{n-1}\right) $ as 
\begin{equation}
\begin{array}{c}
\left( 2k_{n-1}-k\right) \Delta q_{n}-2k_{n-1}\left( k_{n-1}-k\right) \left(
\nabla q_{n}\right) ^{2}-2kh\nabla \overline{q_{n-1}}\nabla q_{n} \\ 
-4h\left( k_{n-1}-k\right) \nabla \overline{q_{n-1}}\nabla q_{n}-2\left(
h\nabla \overline{q_{n-1}}\right) ^{2}+2h\Delta \overline{q_{n-1}} \\ 
+2k\nabla q_{n}\nabla V+4\left( k_{n-1}-k\right) \nabla V\nabla
q_{n}+4\nabla Vh\nabla \overline{q_{n-1}}=2\left( \Delta V+\left( \nabla
V\right) ^{2}\right) , \\ 
q_{n}\mid _{\partial \Omega }=\psi _{n}\left( x\right) .%
\end{array}
\label{4.19}
\end{equation}%
Assuming that the number $h$ is sufficiently small and that $h\overline{k}%
<<1,$ we now ignore those terms in (\ref{4.19}), whose absolute values are $%
O\left( h\right) $ as $h\rightarrow 0$. We also assume in the convergence
analysis that the number $\overline{k}-\underline{k}$ is sufficiently small.
Hence, the number $\left\vert h\nabla \overline{q_{n-1}}\right\vert \leq
\left( \overline{k}-\underline{k}\right) \max_{j}\sup_{x\in \Omega
}\left\vert \nabla q_{j}\left( x\right) \right\vert $ is also small.
However, we do not ignore $\left\vert h\nabla \overline{q_{n-1}}\right\vert
. $ Still, we ignore in (\ref{4.19}) the term $2\left( h\nabla \overline{%
q_{n-1}}\right) ^{2}.$ Hence, we obtain from (\ref{4.19}) for 
\begin{equation}
\begin{array}{c}
k_{n-1}\Delta q_{n}-2k\nabla q_{n}h\nabla \overline{q_{n-1}}+2k\nabla
q_{n}\nabla V+2h\Delta \overline{q_{n-1}} \\ 
=2\left( \Delta V+\left( \nabla V\right) ^{2}\right) -4\nabla Vh\nabla 
\overline{q_{n-1}},\text{ }x\in \Omega ,k\in \left[ k_{n},k_{n-1}\right) ,
\\ 
q_{n}\mid _{\partial \Omega }=\psi _{n}\left( x\right) .%
\end{array}
\label{4.20}
\end{equation}%
Even though the left hand side of equation (\ref{4.20}) depends on $k$, it
changes very little with respect to $k\in \left[ k_{n},k_{n-1}\right) $
since the interval $\left[ k_{n},k_{n-1}\right) $ is small. Still, to
eliminate this $k-$dependence, we integrate both sides of equation (\ref%
{4.20}) with respect to $k\in \left( k_{n},k_{n-1}\right) $ and then divide
both sides of the resulting equation by $h$. We obtain%
\begin{equation}
\begin{array}{c}
\Delta q_{n}-A_{n}h\overline{\nabla q_{n-1}}\nabla q_{n}= \\ 
-A_{n}\nabla q_{n-1}\nabla V_{n-1}+2\left( \Delta V_{n-1}+\left( \nabla
V_{n-1}\right) ^{2}\right) /k_{n-1} \\ 
-4\nabla V_{n-1}h\nabla \overline{q_{n-1}}/k_{n-1}-2h\Delta \overline{q_{n-1}%
}/k_{n-1},\text{ }x\in \Omega , \\ 
q_{n}\mid _{\partial \Omega }=\psi _{n}\left( x\right) ,%
\end{array}
\label{4.21}
\end{equation}%
where $A_{n}=\left( 1+k_{n}/k_{n-1}\right) .$ Hence, 
\begin{equation}
0<A_{n}<2.  \label{4.22}
\end{equation}%
We have replaced in (\ref{4.21}) $V$ with $V_{n-1}$ since we will update
tail functions in our iterative algorithm. In (\ref{4.21}) the term $%
-A_{n}\nabla q_{n-1}\nabla V_{n-1}$ should actually be $-A_{n}\nabla
q_{n}\nabla V_{n-1}.$ We have made this replacement for our convergence
analysis. Indeed, for the exact solution $\nabla q_{n}^{\ast }=\nabla
q_{n-1}^{\ast }+\left( \nabla q_{n}^{\ast }-\nabla q_{n-1}^{\ast }\right) .$
By Lemma 4.3 $\left\vert \nabla q_{n}^{\ast }-\nabla q_{n-1}^{\ast
}\right\vert =O\left( h\right) ,h\rightarrow 0$. The latter, the above
dropped terms, whose absolute values are $O\left( h\right) $ as $%
h\rightarrow 0,$ as well as the approximations (\ref{4.16}) are taken into
account by the function $G_{n}^{\ast }\left( x\right) $ in (\ref{5.11}) and (%
\ref{5.12}) in our convergence analysis.

\begin{algorithm}[Globally convergent algorithm]
\label{Algorithm 1}~

\begin{enumerate}
\item \textrm{\ }

\item \textrm{Set $q_{0}\equiv 0,q_{1,0}=0,\nabla V_{1,0}=\nabla V_{0},$
where the vector $\nabla V_{0}$ is found as in Subsection \ref{subsection
5.3}. }

\item \textrm{For $n = 1$ to $N$, }

\begin{enumerate}
\item \textrm{\ }

\item \textrm{Assume that $q_{n - 1}$ and $\nabla V_{n - 1}$ are found. Set $%
q_{n, 0} = q_{n - 1}$ and $\nabla V_{n, 0} = \nabla V_{n - 1}$. }

\item \textrm{For $i=1$ to $m$ (for an integer $m\geq 1$) }

\begin{enumerate}
\item \textrm{\ }

\item \textrm{Assume that $q_{n,i-1}$ and $V_{n,i-1}$ are found. Find $%
q_{n,i}\in C^{2+\alpha }\left( \overline{\Omega }\right) $ as the solution
of the Dirichlet boundary value problem: }%
\begin{equation}
\begin{array}{c}
\Delta q_{n,i}-A_{n}h\overline{\nabla q_{n-1}}\nabla q_{n,i}= \\ 
-A_{n}\nabla q_{n-1}\nabla V_{n,i-1}+2\left( \Delta V_{n,i-1}+\left( \nabla
V_{n,i-1}\right) ^{2}\right) /k_{n-1} \\ 
-4\nabla V_{n,i-1}h\overline{\nabla q_{n-1}}/k_{n-1}-2h\Delta \overline{%
q_{n-1}}/k_{n-1},\text{ }x\in \Omega , \\ 
q_{n,i}\mid _{\partial \Omega }=\psi _{n}\left( x\right) .%
\end{array}
\label{4.23}
\end{equation}

\item \textrm{Consider the vector function $\nabla v_{n,i}$, 
\begin{equation}
\nabla v_{n,i}\left( x\right) =-\left( hq_{n,i}\left( x\right) +h\overline{%
q_{n-1}}\left( x\right) \right) +\nabla V_{n,i-1}\left( x\right) ,\text{ }%
x\in \Omega .  \label{4.24}
\end{equation}%
Using $\Delta v_{n,i}=\func{div}\left( \nabla v_{n,i}\right) ,$ calculate
the approximation $c_{n,i}\left( x\right) \in C^{\alpha }\left( \overline{%
\Omega }\right) $ for the target coefficient $c\left( x\right) $ as 
\begin{eqnarray}
\beta _{n,i}(x) &=&-\frac{1}{k_{n}^{2}}(\Delta v_{n,i}(x)+(\nabla
v_{n,i}(x))^{2})-1,  \label{4.25} \\
c_{n,i}\left( x\right) &=&\beta _{n,i}(x)+1.  \label{4.26}
\end{eqnarray}%
}

\item \textrm{\label{step update u alg 1} Next, solve the Lippman-Schwinger
equation (\ref{8.8}), where $\widehat{\rho }(y)$ is replaced with $\chi
\left( y\right) (c_{n,i}-1)\left( y\right) .$ We obtain the function $%
u_{n,i}\left( x,\overline{k}\right) .$ Update the first derivatives of the
tail function by 
\begin{equation}
\nabla V_{n,i}\left( x\right) =\frac{\nabla u_{n,i}(x,\overline{k})}{%
u_{n,i}(x,\overline{k})}.  \label{4.250}
\end{equation}%
}
\end{enumerate}
\end{enumerate}

\item \textrm{Set $q_{n}=q_{n,m}$, $c_{n}=c_{n,m}$. }

\item Let $\overline{N}\in \left[ 1,N\right] $ be the optimal number for the
stopping criterion. \textrm{Set the function $c_{\overline{N}}$ as the
computed solution of Problem \ref{Problem 2.1}. }
\end{enumerate}
\end{algorithm}

\textbf{Remarks 5.2}:

\begin{enumerate}
\item The number $\overline{N}\in \left[ 1,N\right] $ should be chosen in
numerical experiments. Our experience with previous works \cite%
{BK1,KFB,KSNF1,TBKF1,TBKF2} indicates that this is possible, also see
section 8. Recall that the number of iteration is often considered as a
regularization parameter in the theory of ill-posed problems, see, e.g. \cite%
{BK1,T}.

\item We solve problems (\ref{4.23}) via the FEM\ using the standard
piecewise linear finite elements. We are doing this, using FreeFem++ \cite{H}%
, which is a very convenient software for the FEM. We note that even though
the Laplacian $\Delta V_{n,i-1}$ is involved in (\ref{4.23}), the Laplacian
is not involved in the variational form of equation (\ref{4.23}). Therefore,
there is no need to calculate $\Delta V_{n,i-1}$ when computing functions $%
q_{n,i}$. Rather, only the gradient $\nabla V_{n,i-1}$ should be calculated.
On the other hand, $\Delta V_{n,i-1}$ should be calculated to find the
function $c_{n,i}(x)$ via (\ref{4.25}), (\ref{4.26}). To calculate $\Delta
V_{n,i-1},$ we use finite differences. The software FreeFem++ automatically
interpolates any function, defined by finite elements to the rectangular
grid and we use this grid to arrange finite differences.

\item Inequality (\ref{3.9}) is valid only if the function $c\left( x\right) 
$ satisfies conditions (\ref{2.0}), (\ref{2.2}) and if Assumption 3.1 holds.
In Theorem 7.1 we impose these conditions on the exact solution $c^{\ast
}\left( x\right) $ of our CIP. However, in the above algorithm, we obtain
functions $c_{n,i}\left( x\right) \in C^{\alpha }\left( \overline{\Omega }%
\right) $ (Theorem 7.1)$,$ which do not necessarily satisfy these
conditions. Nevertheless, we prove in Theorem 7.1 that functions $%
u_{n,i}\left( x,\overline{k}\right) \neq 0,\forall x\in \overline{\Omega }.$
It follows from (\ref{4.250}) that the latter is sufficient for our
algorithm.
\end{enumerate}

\section{Existence and uniqueness of the solution of the Dirichlet boundary
value problem (\protect\ref{4.23})}

\label{sec:6}

In this section we study the question of existence and uniqueness of the
solution $q_{n}\in C^{2+\alpha }\left( \overline{\Omega }\right) $ of the
the Dirichlet boundary value problem (\ref{4.23}). If we would deal with
real valued functions, then existence and uniqueness would follow
immediately from the maximum principle and Schauder theorem \cite{GT,Lad}.
However, complex valued functions cause some additional difficulties, since
the maximum principle does not work in this case. Still, we can get our
desired results using the assumption that the function $\left\vert A_{n}h%
\overline{\nabla q_{n-1}}\right\vert \left( x\right) $ is sufficiently
small, since we assume in Theorem 7.1 that the number $a=\overline{k}-%
\underline{k}=Nh$ is sufficiently small. Keeping in mind the convergence
analysis in the next section, it is convenient to consider here the
Dirichlet boundary value problem (\ref{4.23}) in a more general form, 
\begin{equation}
\begin{array}{c}
\Delta w-\nabla p\nabla w=f\left( x\right) ,\text{ }x\in \Omega , \\ 
w\mid _{\partial \Omega }=\mu \left( x\right) .%
\end{array}
\label{5.1}
\end{equation}

\textbf{Theorem 6.1}. \emph{Assume that in (\ref{5.1}) all functions are
complex valued ones and also that }$p\in C^{1+\alpha }\left( \overline{%
\Omega }\right) ,f\in C^{\alpha }\left( \overline{\Omega }\right) ,\mu \in
C^{2+\alpha }\left( \partial \Omega \right) .$\emph{\ Then there exists a
constant }$C_{1}=C_{1}\left( \Omega \right) >0$\emph{\ depending only on the
domain }$\Omega $ \emph{and a} \emph{sufficiently small number }$\sigma
=\sigma \left( C_{1}\right) \in \left( 0,1\right) $\emph{\ such that if }$%
C_{1}\sigma <1/2$ \emph{and} $\left\vert \nabla p\right\vert _{\alpha }\leq
\sigma $,\emph{\ then there exists unique solution }$w\in C^{2+\alpha
}\left( \overline{\Omega }\right) $\emph{\ of the Dirichlet boundary value
problem (\ref{5.1}) and also }%
\begin{equation}
\left\vert w\right\vert _{2+\alpha }\leq C_{1}\left( \left\vert f\right\vert
_{\alpha }+\left\Vert \mu \right\Vert _{C^{2+\alpha }\left( \partial \Omega
\right) }\right) .  \label{5.2}
\end{equation}

\textbf{Proof}. Below in this paper $C_{1}=C_{1}\left( \Omega \right) >0$
denotes different positive constants depending only on the domain $\Omega .$
Let the complex valued function $v\in C^{2+\alpha }\left( \overline{\Omega }%
\right) .$ Consider the following Dirichlet boundary value problem with
respect to the function $U$:%
\begin{equation}
\begin{array}{c}
\Delta U=\nabla p\nabla v+f\left( x\right) ,\text{ }x\in \Omega , \\ 
U\mid _{\partial \Omega }=\mu \left( x\right) .%
\end{array}
\label{5.3}
\end{equation}%
The Schauder theorem implies that there exists unique solution $U\in
C^{2+\alpha }\left( \overline{\Omega }\right) $ of the problem (\ref{5.3})
and 
\begin{equation}
\left\vert U\right\vert _{2+\alpha }\leq C_{1}\left( \sigma \left\vert
v\right\vert _{2+\alpha }+\left\vert f\right\vert _{\alpha }+\left\Vert \mu
\right\Vert _{C^{2+\alpha }\left( \partial \Omega \right) }\right) .
\label{5.4}
\end{equation}%
Hence, for each fixed pair $f\in C^{\alpha }\left( \overline{\Omega }\right)
,\mu \in C^{2+\alpha }\left( \partial \Omega \right) $, we can define a map
that sends the function $v\in C^{2+\alpha }\left( \overline{\Omega }\right) $
in the solution $U\in C^{2+\alpha }\left( \overline{\Omega }\right) $ of the
problem (\ref{5.3}), say $U=S_{f,\mu }\left( v\right) .$ Hence, $S_{f,\mu
}:C^{2+\alpha }\left( \overline{\Omega }\right) \rightarrow C^{2+\alpha
}\left( \overline{\Omega }\right) .$ Since the operator $S_{f,\mu }$ is
affine and $C_{1}\sigma <1/2,$ then (\ref{5.4}) implies that $S_{f,\mu }$ is
contraction mapping. Let the function $w=S_{f,\mu }\left( w\right) $ be its
unique fixed point. Then the function $w\in C^{2+\alpha }\left( \overline{%
\Omega }\right) $ is the unique solution of the problem (\ref{5.1}). In
addition, by (\ref{5.4}) 
\begin{equation*}
\left\vert w\right\vert _{2+\alpha }\leq C_{1}\left( \sigma \left\vert
w\right\vert _{2+\alpha }+\left\vert f\right\vert _{\alpha }+\left\Vert \mu
\right\Vert _{C^{2+\alpha }\left( \partial \Omega \right) }\right) .
\end{equation*}%
Hence, $\left\vert w\right\vert _{2+\alpha }\leq 2C_{1}\left( \left\vert
f\right\vert _{\alpha }+\left\Vert \mu \right\Vert _{C^{2+\alpha }\left(
\partial \Omega \right) }\right) .$ $\square $

Corollary 6.1 follows immediately from Theorem 6.1, (\ref{4.22}) and (\ref%
{4.23}).

\textbf{Corollary 6.1}. \emph{Consider the Dirichlet boundary value problem (%
\ref{4.23}). Assume that} 
\begin{equation*}
\nabla V_{n,i-1}\in C^{1+\alpha }\left( \overline{\Omega }\right) \text{ and 
}q_{s}\in C^{2+\alpha }\left( \overline{\Omega }\right) ,s=1,...,n-1.
\end{equation*}%
\emph{Suppose that }$\left\vert q_{s}\right\vert _{2+\alpha }\leq Y,$\emph{\
where }$Y=const.>0.$\emph{\ Let }$C_{1}$\emph{\ and }$\sigma $\emph{\ be the
constants of Theorem 6.1 and let the length }$a=\overline{k}-\underline{k}$%
\emph{\ of the interval }$\left[ \underline{k},\overline{k}\right] $\emph{\
be so small that }%
\begin{equation*}
2Ya\leq C_{1}\sigma <1/2.\emph{\ }
\end{equation*}%
\emph{Then there exists unique solution }$q_{n}\in C^{2+\alpha }\left( 
\overline{\Omega }\right) $\emph{\ of the problem (\ref{4.23}) and}%
\begin{equation*}
\begin{array}{c}
\emph{\ }\left\vert q_{n}\right\vert _{2+\alpha }\leq C_{1}\left\vert
A_{n}\nabla q_{n-1}\nabla V_{n,i-1}\right\vert _{\alpha }+\left\vert 2\left(
\Delta V_{n,i-1}+\left( \nabla V_{n,i-1}\right) ^{2}\right)
/k_{n-1}\right\vert _{\alpha } \\ 
+C_{1}\left\vert 4\nabla V_{n,i-1}h\overline{\nabla q_{n-1}}%
/k_{n-1}\right\vert _{\alpha }+C_{1}\left\vert 2h\Delta \overline{q_{n-1}}%
/k_{n-1}\right\vert _{\alpha }+C_{1}\left\Vert \psi _{n}\right\Vert
_{C^{2+\alpha }\left( \partial \Omega \right) }.%
\end{array}%
\end{equation*}

\section{Global convergence}

\label{sec:7}

Let $\delta >0$ be the level of the error in the boundary data $\psi
_{n}\left( x\right) $ in (\ref{4.16}). We introduce the error parameter $%
\eta ,$ 
\begin{equation}
\eta =h+\delta .  \label{5.7}
\end{equation}%
It is natural to assume that 
\begin{equation}
\left\Vert \psi _{n}-\psi _{n}^{\ast }\right\Vert _{C^{2+\alpha }\left(
\partial \Omega \right) }\leq \eta .  \label{5.8}
\end{equation}%
We now introduce some natural assumptions about the exact coefficient $%
c^{\ast }\left( x\right) $ and functions associated with it. Let $%
B_{1}^{\ast }=B_{1}^{\ast }\left( \Omega ,c^{\ast }\right) >0$ be the number
of Lemma 3.1, which corresponds to $c^{\ast }.$ We assume that the number $%
\underline{k}>0$ is so large that 
\begin{equation}
\frac{B_{1}^{\ast }}{\underline{k}}<\frac{1}{2}.  \label{200}
\end{equation}%
Let $A^{\ast }\left( x\right) >0$ be the function $A\left( x\right) $ in (%
\ref{3.8}) which corresponds to $c^{\ast }.$ Denote%
\begin{equation}
A_{\min }^{\ast }=\min_{\overline{\Omega }}A^{\ast }\left( x\right) ,\text{ }%
D^{\ast }=\min \left[ \left( A_{\min }^{\ast }\right) ^{2},\left( A_{\min
}^{\ast }\right) ^{4}\right] .  \label{7.8}
\end{equation}%
Hence, (\ref{3.8}), (\ref{200}) and (\ref{7.8}) imply that 
\begin{equation}
\min_{\overline{\Omega }}\left\vert u^{\ast }\left( x,k\right) \right\vert
\geq \frac{A_{\min }^{\ast }}{2},\forall k\geq \underline{k}.  \label{201}
\end{equation}

Following (\ref{def q}), let $q^{\ast }\left( x,k\right) =\partial
_{k}u^{\ast }\left( x,k\right) /u^{\ast }\left( x,k\right) .$ All
approximations for the function $q^{\ast }\left( x,k\right) $ and associated
functions with the accuracy $O\left( h\right) $ as $h\rightarrow 0,$ which
are used in this section below, can be justified by Lemma 4.3. Denote $%
q_{n}^{\ast }\left( x\right) =q^{\ast }\left( x,k_{n}\right) .$ For $k\in %
\left[ k_{n},k_{n-1}\right) $ we obtain $q^{\ast }\left( x,k\right)
=q_{n}^{\ast }\left( x\right) +O\left( h\right) $ as $h\rightarrow 0.$ Set $%
q_{0}^{\ast }\left( x\right) \equiv 0.$ Using (\ref{4.24}), define the
gradient $\nabla v_{n}^{\ast }$ as 
\begin{equation}
\nabla v_{n}^{\ast }\left( x\right) =-h\nabla {q_{n}^{\ast }}\left( x\right)
-h\nabla \overline{q_{n-1}^{\ast }}\left( x\right) +\nabla V^{\ast }\left(
x\right) ,\text{ }x\in \Omega .  \label{5.9}
\end{equation}%
Since $\Delta v_{n}^{\ast }=\func{div}\left( \nabla v_{n}^{\ast }\right) ,$
then by (\ref{4.1})%
\begin{equation}
c^{\ast }\left( x\right) =-\frac{1}{k_{n}^{2}}\left( \Delta v_{n}^{\ast
}+\left( \nabla v_{n}^{\ast }\right) ^{2}\right) +F_{n}^{\ast }\left(
x\right) .  \label{5.10}
\end{equation}%
While equation (\ref{4.6}) is precise, we have obtained equation (\ref{4.19}%
) using some approximations whose error is $O\left( h\right) $ as $%
h\rightarrow 0.$ This justifies the presence of the term $G_{n}^{\ast
}\left( x\right) $ in the following analog of the Dirichlet boundary value
problem (\ref{4.19}):%
\begin{equation}
\begin{array}{c}
\Delta q_{n}^{\ast }-A_{n}h\overline{\nabla q_{n-1}^{\ast }}\nabla
q_{n}^{\ast }= \\ 
-A_{n}\nabla q_{n-1}^{\ast }\nabla V^{\ast }+2\left( \Delta V^{\ast }+\left(
\nabla V^{\ast }\right) ^{2}\right) /k_{n-1} \\ 
-4\nabla V^{\ast }h\nabla \overline{q_{n-1}^{\ast }}/k_{n-1}-2h\Delta 
\overline{q_{n-1}^{\ast }}/k_{n-1}+G_{n}^{\ast }\left( x\right) ,\text{ }%
x\in \Omega , \\ 
q_{n}^{\ast }\mid _{\partial \Omega }=\psi _{n}^{\ast }\left( x\right) .%
\end{array}
\label{5.11}
\end{equation}%
In (\ref{5.10}), (\ref{5.11}) $F_{n}^{\ast }\left( x\right) $ and $%
G_{n}^{\ast }\left( x\right) $ are error functions, which can be estimated as%
\begin{equation}
\left\vert F_{n}^{\ast }\right\vert _{\alpha }\leq M\eta ,\left\vert
G_{n}^{\ast }\right\vert _{\alpha }\leq M\eta ,  \label{5.12}
\end{equation}%
where $M>0$ is a constant. We also assume that 
\begin{equation}
\left\vert \nabla V^{\ast }\right\vert _{\alpha },\left\vert \Delta V^{\ast
}\right\vert _{\alpha },\left\vert q_{n}^{\ast }\right\vert _{2+\alpha
},\left\vert \nabla v_{n}^{\ast }\right\vert _{\alpha },\left\vert \Delta
v_{n}^{\ast }\right\vert _{\alpha },\left\vert u^{\ast }\left( x,\overline{k}%
\right) \right\vert _{2+\alpha }\leq M.  \label{5.13}
\end{equation}

Denote $C_{2}=C_{2}\left( \Omega \right) =\max \left( C\left( \Omega \right)
,C_{1}\left( \Omega \right) \right) >0,$ where $C\left( \Omega \right) $ and 
$C_{1}\left( \Omega \right) $ are constants of Theorem 5.1 and Corollary 6.1
respectively. For brevity and also to emphasize the main idea of the proof,
we formulate and prove Theorem 7.1 only for the case when inner iterations
are absent, i.e. for the case $m=0$. The case $m\geq 1$ is a little bit more
technical and is, therefore, more space consuming, while the idea is still
the same. Thus, any function $f_{nj}$ in above formulae should be $f_{n}$
below in this section.

\textbf{Theorem 7.1} (global convergence). \emph{Let the exact coefficient }$%
c^{\ast }\left( x\right) =\beta ^{\ast }\left( x\right) +1$\emph{\ satisfies
conditions (\ref{2.0}), (\ref{2.2}) and let Assumption 3.1 holds for }$%
c^{\ast }\left( x\right) $\emph{. Assume that the first approximation }$%
V_{0}\left( x\right) $\emph{\ for the tail function is constructed as in
subsection 5.3. Let numbers }$\overline{k}>\underline{k}>1$\emph{\ and let} 
\emph{(\ref{5.7})-(\ref{200}) hold true. Let the number }$Z^{\ast }=Z^{\ast
}\left( \beta ^{\ast },\overline{k},\chi ,\Omega _{1},\Omega \right) >0$%
\emph{\ and a sufficiently small number }$\theta ^{\ast }=\theta ^{\ast
}\left( \beta ^{\ast },\overline{k},\chi ,\Omega _{1},\Omega \right) \in
\left( 0,1\right) $\emph{\ be the constants of Theorem 4.1, which depend
only on listed parameters. Assume that the number }$M$\emph{\ in (\ref{5.12}%
) and (\ref{5.13}) is so large that }

\begin{equation}
M>\max \left( 4,Z^{\ast },\frac{256}{D^{\ast }},28C_{2},C_{2}\overline{k}%
\right) .  \label{5.131}
\end{equation}%
\emph{\ Let the number }$a=\overline{k}-\underline{k}$\emph{\ be so small
that }%
\begin{equation}
4Ma<C_{1}\sigma <1/2,  \label{5.16}
\end{equation}%
\emph{where }$C_{1}$\emph{\ and }$\sigma $\emph{\ are numbers of Theorem 6.1.%
} \emph{Let $N\geq 2$ and the level of the error }$\eta $\emph{\ in (\ref%
{5.7}) be such that }%
\begin{equation}
\eta \in \left( 0,\eta _{0}\right) ,\text{ }\eta _{0}=\frac{\theta ^{\ast }}{%
M^{20N-12}}.  \label{5.130}
\end{equation}%
\emph{Assume that the number }$\theta ^{\ast }$\emph{\ is so small that }%
\begin{equation}
\theta ^{\ast }<\frac{A_{\min }^{\ast }}{4}.  \label{5.132}
\end{equation}%
\emph{Then for }$n=1,2,...,N$\emph{\ reconstructed functions }$c_{n}\in
C^{\alpha }\left( \overline{\Omega }\right) $\emph{\ and also }$\min_{%
\overline{\Omega }}\left\vert u_{n}\left( x,\overline{k}\right) \right\vert
\geq A_{\min }^{\ast }/4.$ \emph{In addition, the following accuracy
estimate is valid }%
\begin{equation}
\left\vert c_{n}-c^{\ast }\right\vert _{\alpha }\leq M^{10N-6}\eta <\sqrt{%
\eta }.  \label{600}
\end{equation}

\textbf{Remark 7.1. }Thus, Theorem 7.1 claims that our iteratively found
functions $c_{n}$\ are located in a sufficiently small neighborhood of the
exact solution $c^{\ast },$\ as long as $n\in \left[ 1,N\right] $ and the
error parameter $\eta $ is sufficiently small. This is achieved without any
advanced knowledge of a small neighborhood of the exact solution $c^{\ast }$%
. Hence, Theorem 7.1 implies the global convergence of our algorithm, see
Introduction. On the other hand, this is achieved within the framework of
the approximation of subsection 5.3. Hence, to be more precise, this is the 
\emph{approximate global convergence }property as defined in \cite{BK1,KSNF1}%
. It can be seen from the proof of this theorem that the approximations
approximations (\ref{4.8}), (\ref{102}) and (\ref{103}) are not used on
follow up iterations with $n=1,...,N$. Recall that the number of iterations (%
$N$\ in our case) can be considered sometimes as a regularization parameter
in the theory of ill-posed problems \cite{BK1,T}. Also, see Remarks 5.1.

\textbf{Proof}. Denote%
\begin{equation}
\begin{array}{c}
D^{\gamma }\widetilde{V}_{n}=D^{\gamma }V_{n}-D^{\gamma }V^{\ast },%
\widetilde{q}_{n}=q_{n}-q_{n}^{\ast },\widetilde{v}_{n}=v_{n}-v_{n}^{\ast },
\\ 
\widetilde{u}_{n}\left( x\right) =u_{n}\left( x,\overline{k}\right) -u^{\ast
}\left( x,\overline{k}\right) ,\widetilde{c}_{n}=c_{n}-c^{\ast },\widetilde{%
\psi }_{n}=\psi _{n}-\psi _{n}^{\ast }.%
\end{array}
\label{5.14}
\end{equation}%
Here $\gamma =\left( \gamma _{1},\gamma _{2},\gamma _{3}\right) $ is multi
index with non-negative integer components and $\left\vert \gamma
\right\vert =\gamma _{1}+\gamma _{2}+\gamma _{3}.$ Using (\ref{107}), (\ref%
{5.7}), (\ref{5.8}) and (\ref{5.131}), we obtain 
\begin{equation}
\left\vert \nabla \widetilde{V}_{0}\right\vert _{\alpha },\left\vert \Delta 
\widetilde{V}_{0}\right\vert _{\alpha }\leq C_{2}\overline{k}\eta \leq M\eta
.  \label{5.140}
\end{equation}%
Hence, (\ref{5.13}), (\ref{5.130}) and (\ref{5.140}) imply that 
\begin{equation}
\left\vert \nabla V_{0}\right\vert _{\alpha }=\left\vert \nabla \widetilde{V}%
_{0}+\nabla V^{\ast }\right\vert _{\alpha }\leq M\eta +M\leq 2M,\text{ }%
\left\vert \Delta V_{0}\right\vert _{\alpha }\leq 2M.  \label{5.1400}
\end{equation}%
Subtract equation (\ref{5.11}) from equation (\ref{4.21}). Also, subtract
the boundary condition in (\ref{5.11}) from the boundary condition in (\ref%
{4.21}). We obtain%
\begin{equation}
\begin{array}{c}
\Delta \widetilde{q}_{n}-A_{n}h\nabla \overline{q_{n-1}}\nabla \widetilde{q}%
_{n}=\widetilde{Q}_{n}, \\ 
\widetilde{Q}_{n}=A_{n}h\nabla \overline{\widetilde{q}_{n-1}}\nabla
q_{n}^{\ast }-A_{n}\left( \nabla \widetilde{q}_{n-1}\nabla V_{n-1}+\nabla
q_{n-1}^{\ast }\nabla \widetilde{V}_{n-1}\right) \\ 
+2\left( \Delta \widetilde{V}_{n-1}+\nabla \widetilde{V}_{n-1}\left( \nabla
V_{n-1}+\nabla V^{\ast }\right) \right) /k_{n-1} \\ 
-4\nabla V_{n-1}h\nabla \overline{\widetilde{q}_{n-1}}/k_{n-1}-4\nabla 
\widetilde{V}_{n-1}h\nabla \overline{q_{n-1}^{\ast }}/k_{n-1}-2h\Delta 
\overline{\widetilde{q}_{n-1}}/k_{n-1}-G_{n}^{\ast }, \\ 
\widetilde{q}_{n}\mid _{\partial \Omega }=\widetilde{\psi }_{n}.%
\end{array}
\label{5.15}
\end{equation}%
Let $n\geq 2$ $\ $and let an integer $p_{n-1}\in \left[ 1,10\left(
N-1\right) -6\right] .$ Assume that 
\begin{equation}
\left\vert \nabla \widetilde{V}_{n-1}\right\vert _{\alpha },\left\vert
\Delta \widetilde{V}_{n-1}\right\vert _{\alpha },\left\vert \widetilde{q}%
_{s}\right\vert _{2+\alpha }\leq M^{p_{n-1}}\eta \leq M,  \label{5.18}
\end{equation}%
where $s=1,...,n-1.$ Note that while the left inequality (\ref{5.18}) is our
assumption, the right inequality (\ref{5.18}) follows from (\ref{5.130}).
Similarly with (\ref{5.1400}) we obtain from (\ref{5.18}) 
\begin{equation}
\left\vert \nabla V_{n-1}\right\vert _{\alpha },\text{ }\left\vert \Delta
V_{n-1}\right\vert _{\alpha },\left\vert q_{s}\right\vert _{2+\alpha }\leq
2M.  \label{5.22}
\end{equation}%
It follows from (\ref{4.22}) and (\ref{5.22}) that 
\begin{equation}
\left\vert A_{n}h\overline{\nabla q_{n-1}}\right\vert _{1+\alpha }\leq 4Ma.
\label{5.220}
\end{equation}%
Hence, Corollary 6.1, (\ref{5.8}), (\ref{5.16}), (\ref{5.15}), (\ref{5.22})
and (\ref{5.220}) imply that%
\begin{equation}
\left\vert \widetilde{q}_{n}\right\vert _{2+\alpha }\leq C_{2}\left\vert 
\widetilde{Q}_{n}\right\vert _{\alpha }+C_{2}\eta .  \label{5.17}
\end{equation}%
We now want to find the number $p_{n}.$ First, using (\ref{5.13}), (\ref%
{5.16}), (\ref{5.15}), (\ref{5.18}) and (\ref{5.22}) and also recalling that
by (\ref{4.22}) $A_{n}<2$, we estimate $\left\vert \widetilde{Q}%
_{n}\right\vert _{\alpha },$%
\begin{equation*}
\begin{array}{c}
\left\vert \widetilde{Q}_{n}\right\vert _{\alpha }\leq \left( 2Ma\right)
M^{p_{n-1}}\eta +6MM^{p_{n-1}}\eta +8MM^{p_{n-1}}\eta +8MM^{p_{n-1}}\eta \\ 
+\left( 8Ma\right) M^{p_{n-1}}\eta +\left( 6Ma\right) M^{p_{n-1}}\eta +M\eta
\\ 
\leq 2M^{p_{n-1}}\eta +24MM^{p_{n-1}}\eta +M\eta \leq 27MM^{p_{n-1}}\eta .%
\end{array}%
\end{equation*}%
Hence, using (\ref{5.8}) and (\ref{5.17}), we obtain $\left\vert \widetilde{q%
}_{n}\right\vert _{2+\alpha }\leq 27C_{2}MM^{p_{n-1}}\eta +C_{2}\eta \leq
28C_{2}MM^{p_{n-1}}\eta .$ Since by (\ref{5.131}) $M>28C_{2},$ then 
\begin{equation}
\left\vert \widetilde{q}_{n}\right\vert _{2+\alpha }\leq M^{p_{n-1}+2}\eta .
\label{5.23}
\end{equation}%
Hence, 
\begin{equation}
\left\vert q_{n}\right\vert _{2+\alpha }\leq \left\vert \widetilde{q}%
_{n}\right\vert _{2+\alpha }+\left\vert q_{n}^{\ast }\right\vert _{2+\alpha
}\leq 2M.  \label{5.24}
\end{equation}

Subtracting (\ref{5.9}) from (\ref{4.24}) and using (\ref{5.14}), we obtain%
\begin{equation*}
\nabla \widetilde{v}_{n}=-h\nabla \widetilde{q}_{n}-h\overline{\nabla 
\widetilde{q}_{n-1}}+\nabla \widetilde{V}_{n-1}.
\end{equation*}%
Hence, using (\ref{5.7}), (\ref{5.18}) and (\ref{5.23}), we obtain 
\begin{equation}
\left\vert \nabla \widetilde{v}_{n}\right\vert _{\alpha },\left\vert \Delta 
\widetilde{v}_{n}\right\vert _{\alpha }\leq M^{p_{n-1}+2}\eta ^{2}+\left(
Ma\right) M^{p_{n-1}}\eta +M^{p_{n-1}}\eta \leq 2M^{p_{n-1}}\eta \leq
M^{p_{n-1}+1}\eta .  \label{5.25}
\end{equation}%
Hence, using (\ref{5.13}) and (\ref{5.25}), we obtain 
\begin{equation}
\left\vert \nabla v_{n}\right\vert _{\alpha }\leq \left\vert \nabla 
\widetilde{v}_{n}\right\vert _{\alpha }+\left\vert \nabla v_{n}^{\ast
}\right\vert _{\alpha }\leq 2M.  \label{5.26}
\end{equation}%
Next, by (\ref{4.25}), (\ref{4.26}), (\ref{5.10}) and (\ref{5.14}) 
\begin{equation}
\widetilde{c}_{n}=-\frac{1}{k_{n}^{2}}\left( \Delta \widetilde{v}_{n}+\left(
\nabla v_{n}+\nabla v_{n}^{\ast }\right) \nabla \widetilde{v}_{n}\right)
-F_{n}^{\ast }.  \label{5.27}
\end{equation}%
In particular, since the right hand side of (\ref{5.27}) belongs to $%
C^{\alpha }\left( \overline{\Omega }\right) ,$ then the function $\widetilde{%
c}_{n}\in C^{\alpha }\left( \overline{\Omega }\right) .$ Recalling that $%
k_{n}^{2}\geq \underline{k}^{2}>1$ and using (\ref{5.12}), (\ref{5.24})-(\ref%
{5.27}), we obtain%
\begin{equation}
\left\vert \widetilde{c}_{n}\right\vert _{\alpha }\leq \left( 3M+1\right)
M^{p_{n-1}+1}\eta +M\eta \leq 4MM^{p_{n-1}+1}\eta \leq M^{p_{n-1}+2}\eta .
\label{5.28}
\end{equation}

We now estimate $\left\vert \widetilde{u}_{n}\right\vert _{2+\alpha }.$ It
follows from (\ref{8.9}), (\ref{5.130}) and (\ref{5.28}) that%
\begin{equation}
\left\vert \widetilde{u}_{n}\right\vert _{2+\alpha }=\left\vert u_{n}\left(
x,\overline{k}\right) -u^{\ast }\left( x,\overline{k}\right) \right\vert
_{2+\alpha }\leq Z^{\ast }M^{p_{n-1}+2}\eta \leq M^{p_{n-1}+3}\eta .
\label{5.29}
\end{equation}%
Hence, similarly with (\ref{5.26})%
\begin{equation}
\left\vert u_{n}\left( x,\overline{k}\right) \right\vert _{2+\alpha }\leq 2M.
\label{5.30}
\end{equation}%
Now we can estimate $\left\vert u_{n}\left( x,\overline{k}\right)
\right\vert $ from the below. Using (\ref{201}), (\ref{5.130}), (\ref{5.132}%
) and (\ref{5.29}), we obtain 
\begin{equation}
\left\vert u_{n}\left( x,\overline{k}\right) \right\vert \geq \left\vert
u^{\ast }\left( x,\overline{k}\right) \right\vert -\left\vert \widetilde{u}%
_{n}\left( x\right) \right\vert \geq \frac{A_{\min }^{\ast }}{2}%
-M^{p_{n-1}+3}\eta \geq \frac{A_{\min }^{\ast }}{4},\text{ }\forall x\in 
\overline{\Omega }.  \label{5.32}
\end{equation}

We now are ready to estimate H\"{o}lder norms of $\nabla \widetilde{V}%
_{n},\Delta \widetilde{V}_{n},\nabla V_{n},\Delta V_{n}.$ It is obvious that
for any two complex valued functions $f_{1},f_{2}\in C^{\alpha }\left( 
\overline{\Omega }\right) $ such that $f_{2}\left( x\right) \neq 0$ in $%
\overline{\Omega }$ 
\begin{equation}
\left\vert \frac{f_{1}}{f_{2}}\right\vert _{\alpha }\leq \frac{\left\vert
f_{1}\right\vert _{\alpha }\left\vert f_{2}\right\vert _{\alpha }}{%
\left\vert f_{2}\right\vert _{\min }^{2}},\quad \text{where }\>\left\vert
f_{2}\right\vert _{\min }=\min_{\overline{\Omega }}\left\vert
f_{2}\right\vert .  \label{5.33}
\end{equation}%
We have%
\begin{equation}
\nabla \widetilde{V}_{n}=\frac{\nabla u_{n}\left( x,\overline{k}\right) }{%
u_{n}\left( x,\overline{k}\right) }-\frac{\nabla u^{\ast }\left( x,\overline{%
k}\right) }{u^{\ast }\left( x,\overline{k}\right) }=  \label{5.34}
\end{equation}%
\begin{equation*}
\frac{u^{\ast }\left( x,\overline{k}\right) \nabla \left( u_{n}\left( x,%
\overline{k}\right) -u^{\ast }\left( x,\overline{k}\right) \right) +\left(
u^{\ast }\left( x,\overline{k}\right) -u_{n}\left( x,\overline{k}\right)
\right) \nabla u^{\ast }\left( x,\overline{k}\right) }{u_{n}\left( x,%
\overline{k}\right) u^{\ast }\left( x,\overline{k}\right) }
\end{equation*}%
Hence, using (\ref{7.8}), (\ref{201}), (\ref{5.13}), (\ref{5.131}), (\ref%
{5.29}), (\ref{5.32}) and (\ref{5.33}), we obtain%
\begin{equation}
\left\vert \nabla \widetilde{V}_{n}\right\vert _{\alpha }\leq \frac{256}{%
D^{\ast }}M^{p_{n-1}+5}\eta \leq M^{p_{n-1}+6}\eta .  \label{5.35}
\end{equation}%
Next,%
\begin{equation}
\Delta \widetilde{V}_{n}=\left( \frac{\Delta u_{n}}{u_{n}}-\frac{\Delta
u^{\ast }}{u^{\ast }}\right) \left( x,\overline{k}\right) -\left( \frac{%
\nabla u_{n}}{u_{n}}-\frac{\nabla u^{\ast }}{u^{\ast }}\right) \left( x,%
\overline{k}\right) \cdot \left( \frac{\nabla u_{n}}{u_{n}}+\frac{\nabla
u^{\ast }}{u^{\ast }}\right) \left( x,\overline{k}\right) .  \label{5.36}
\end{equation}%
We now estimate each term in (\ref{5.36}). Using the similarity with (\ref%
{5.34}) as well as (\ref{5.35}), we obtain 
\begin{equation}
\left\vert \frac{\Delta u_{n}\left( x,\overline{k}\right) }{u_{n}\left( x,%
\overline{k}\right) }-\frac{\Delta u^{\ast }\left( x,\overline{k}\right) }{%
u^{\ast }\left( x,\overline{k}\right) }\right\vert _{\alpha }\leq
M^{p_{n-1}+6}\eta .  \label{5.37}
\end{equation}%
Next, using (\ref{7.8}), (\ref{5.13}), (\ref{5.30})-(\ref{5.33}), we obtain%
\begin{equation}
\left\vert \frac{\nabla u_{n}\left( x,\overline{k}\right) }{u_{n}\left( x,%
\overline{k}\right) }+\frac{\nabla u^{\ast }\left( x,\overline{k}\right) }{%
u^{\ast }\left( x,\overline{k}\right) }\right\vert _{\alpha }\leq \frac{66}{%
D^{\ast }}M^{2}.  \label{5.370}
\end{equation}%
Hence, using (\ref{5.131}) and (\ref{5.36})-(\ref{5.370}), we obtain%
\begin{equation}
\left\vert \Delta \widetilde{V}_{n}\right\vert _{\alpha }\leq \left( \frac{66%
}{D^{\ast }}M^{2}+1\right) M^{p_{n-1}+6}\eta \leq \left( M^{3}+1\right)
M^{p_{n-1}+6}\eta \leq M^{p_{n-1}+10}\eta .  \label{5.38}
\end{equation}%
Similarly with the above, we derive from (\ref{5.35}) and (\ref{5.38}) that $%
\left\vert \nabla V_{n}\right\vert _{\alpha },\left\vert \Delta
V_{n}\right\vert _{\alpha }\leq 2M.$ Summarizing, assuming the validity of
estimates (\ref{5.18}), we have established the following estimates:%
\begin{equation}
\left\vert \nabla \widetilde{V}_{n}\right\vert _{\alpha },\left\vert \Delta 
\widetilde{V}_{n}\right\vert _{\alpha },\left\vert \widetilde{q}%
_{n}\right\vert _{2+\alpha },\left\vert \nabla \widetilde{v}_{n}\right\vert
_{\alpha },\left\vert \Delta \widetilde{v}_{n}\right\vert _{\alpha
},\left\vert \widetilde{c}_{n}\right\vert _{\alpha }\leq M^{p_{n-1}+10}\eta ,
\label{5.39}
\end{equation}%
\begin{equation*}
\left\vert \nabla V_{n}\right\vert _{\alpha },\left\vert \Delta
V_{n}\right\vert _{\alpha },\left\vert q_{n}\right\vert _{2+\alpha
},\left\vert \nabla v_{n}\right\vert _{1+\alpha },\left\vert \Delta
v_{n}\right\vert _{\alpha }\leq 2M.
\end{equation*}%
Hence, it follows from (\ref{5.18}) and (\ref{5.39}) that $p_{n}=p_{n-1}+10.$
Hence, $p_{n}=p_{1}+10\left( n-1\right) .$ We now need to find $p_{1}.$ Let $%
n=1$. Then (\ref{5.12}), (\ref{5.131}), (\ref{5.140})-(\ref{5.15}) and (\ref%
{5.17}) imply that 
\begin{equation}
\left\vert \widetilde{q}_{1}\right\vert _{2+\alpha }\leq C_{2}\left\vert 
\widetilde{Q}_{1}\right\vert _{\alpha }+C_{2}\eta \leq 2C_{2}\left(
1+3M\right) M\eta +M\eta \leq M^{4}\eta .  \label{5.40}
\end{equation}%
Hence, (\ref{5.140}), (\ref{5.18}) and (\ref{5.40}) imply that $p_{1}=4.$
Hence, $p_{n}=10n-6.$ Thus, estimates (\ref{5.39}) are valid for $%
p_{n-1}+10=10n-6.$ The target estimate (\ref{600}) of this theorem is
equivalent with the estimate for $\left\vert \widetilde{c}_{n}\right\vert
_{\alpha }$ in (\ref{5.39}). $\ \ \ \square $

\section{Numerical study}

\label{sec:8}

In this section, we present numerical \ results for our method. It is well
known that numerical implementations of algorithms quite often deviate
somewhat from the theory. In other words, discrepancies between the theory
and its numerical implementation occur quite often, including this paper.
Since in our target application to imaging of explosives (section 1) only
backscattering data are measured \cite{KSNF1,TBKF1,TBKF2}, we slightly
modify our Algorithm \ref{Algorithm 1} to work with these data. The
computations were performed using the above mentioned (Remark 5.2) software
FreeFem++ \cite{H}.

\subsection{Numerical solution of the forward problem}

\label{sec:8.1}

To generate our data (\ref{2.6}) for the inverse problem, we need to solve
the forward problem (\ref{2.3})-(\ref{2.5}). We solve it via the FEM. Let $%
A>0$ be a number. We set $\Omega =\left( -A,A\right) ^{3}.$ Taking $\Omega $
as a cube is convenient for our planned future work with experimental data,
as the above mentioned past experience of working with time resolved
experimental data demonstrates \cite{TBKF1,TBKF2}. Let the number $A_{1}>A.$
Since it is impossible to numerically solve equation (\ref{2.3}) in the
entire space $\mathbb{R}^{3},$ we \textquotedblleft truncate" this space and
solve this equation in the cube $G=\left( -A_{1},A_{1}\right) ^{3}.$ Hence, $%
\Omega \subset G$ and $\partial \Omega \cap \partial G=\varnothing .$
Consider different parts of the boundaries $\partial G$ and $\partial \Omega
,$%
\begin{equation*}
\begin{array}{c}
\partial G=\partial _{1}G\cup \partial _{2}G\cup \partial _{3}G,\partial
\Omega =\partial _{1}\Omega \cup \partial _{2}\Omega \cup \partial
_{3}\Omega , \\ 
\partial _{1}G=\left\{ x_{1},x_{2}\in \left( -A_{1},A_{1}\right)
,x_{3}=-A_{1}\right\} ,\partial _{1}\Omega =\left\{ x_{1},x_{2}\in \left(
-A,A\right) ,x_{3}=-A\right\} , \\ 
\partial _{2}G=\left\{ x_{1},x_{2}\in \left( -A_{1},A_{1}\right)
,x_{3}=A_{1}\right\} ,\partial _{2}\Omega =\left\{ x_{1},x_{2}\in \left(
-A,A\right) ,x_{3}=A\right\} , \\ 
\partial _{3}G=\partial G\diagdown \left( \partial _{1}G\cup \partial
_{2}G\right) ,\partial _{3}\Omega =\partial \Omega \diagdown \left( \partial
_{1}\Omega \cup \partial _{2}\Omega \right) .%
\end{array}%
\end{equation*}

To generate the data for the inverse problem, we solve the following forward
problems in the cube $G$ for $k=k_{0},...k_{N}:$ 
\begin{equation}
\begin{array}{c}
\Delta u+k^{2}c\left( x\right) u=0\text{ in }G, \\ 
\partial _{n}u+iku=0,x\in \partial _{1}G\cup \partial _{2}G, \\ 
\partial _{n}u=0,x\in \partial _{3}G, \\ 
u=\exp \left( -ikx_{3}\right) +u_{sc}.%
\end{array}
\label{8.1}
\end{equation}%
Recall that $c\left( x\right) =1$ outside of\ the domain $\Omega .$ The
second line in (\ref{8.1}) is the absorbing boundary condition. The
condition in the third line of (\ref{8.1}) can be interpreted as follows:
the vertical boundary $\partial _{3}G$ is so far from inhomogeneities, which
are located in the cube $\Omega ,$ that they do not affect the incident
plane wave $\exp \left( -ikx_{3}\right) $ for $x\in \partial _{3}G.$

\subsection{Backscattering data}

\label{sec:8.2}

Our numerical examples are only for the case of the backscattering data. Let
the number $s\in \left( 0,A\right] .$ Denote 
\begin{equation}
P_{-s}=\left\{ x\in \Omega :x_{3}=-s\right\} ,s>0.  \label{8.0}
\end{equation}%
Hence, the set $P_{-A}=\partial _{1}\Omega $ is the bottom boundary of $%
\Omega .$ We assume that the backscattering data $g\left( x,k\right) $ are
measured on $P_{-A}$. Similarly with \cite{TBKF1,TBKF2} we complement the
data at $P_{-R}$ with the data for the case $c\left( x\right) \equiv 1.$ In
other words, we use the function $\widetilde{g}\left( x,k\right) $ instead
of the function $g\left( x,k\right) ,$ where%
\begin{equation}
\widetilde{g}\left( x,k\right) =\left\{ 
\begin{array}{c}
g\left( x,k\right) ,x\in P_{-A}, \\ 
\exp \left( -ikx_{3}\right) ,x\in \partial \Omega \diagdown P_{-A}.%
\end{array}%
\right.  \label{8.2}
\end{equation}%
Formula (\ref{8.2}) can be intuitively justified in the case when
explosive-like targets of interest are located far from the part of the
boundary $\partial \Omega \diagdown P_{-A}$ of the domain $\Omega .$ To be
in an agreement with Theorem 7.1, we assume that the function $g^{\ast
}\left( x,k\right) \mid _{\partial \Omega \diagdown P_{-A}}=u^{\ast }\left(
x,k\right) \mid _{\partial \Omega \diagdown P_{-A}},$ generated by the exact
coefficient $c^{\ast }\left( x\right) ,$ is close to the function $\exp
\left( -ikx_{3}\right) \mid _{\partial \Omega \diagdown P_{-A}}.$

\subsection{Some details of the numerical implementation}

\label{sec:8.3}

In this subsection we describe some details of the numerical implementation
of Algorithm 5.1.

\subsubsection{Computations of tail functions}

\label{sec:8.3.1}

As it is clear from (\ref{4.23}) and item 2 of Remarks 5.2, we use only the
gradient of each tail function. Recall that by (\ref{106}) the first tail
function $V_{0}\left( x\right) =-i\overline{k}\tau \left( x\right) ,$ where
the function $\tau \left( x\right) $ is the solution of the boundary value
problem (\ref{1060}). Hence, to avoid the noise linked with the
differentiation of $V_{0}\left( x\right) ,$ we have numerically solved the
following problem to calculate the gradient $\nabla V_{0}$: 
\begin{equation}
\begin{array}{c}
\Delta \left( \nabla V_{0}\right) =0,x\in \Omega , \\ 
\nabla V_{0}\left( x\right) =\nabla u\left( x,\overline{k}\right) /u\left( x,%
\overline{k}\right) ,x\in \partial \Omega .%
\end{array}
\label{8.300}
\end{equation}%
Indeed, since by (\ref{106}) and (\ref{1060}) the function $V_{0}$ satisfies
the Laplace equation, then its derivatives also satisfy this equation. The
next question is on how to obtain the boundary data for $\nabla u\left( x,%
\overline{k}\right) .$ There are two ways of doing this. The first way is to
solve equation (\ref{8.1}) in the domain $G^{\prime }=G\diagdown \Omega $
for $k:=\overline{k}$ with the same boundary conditions on $\partial G$ as
in (\ref{8.1}) and with the boundary condition (\ref{8.2}) on $\partial
\Omega .$ In doing so, one should assume that there exists unique solution
of this boundary value problem. In the case when $\partial \Omega \in
C^{2+\alpha },$ which was considered in sections 2-7, one can use (\ref{2.7}%
), (\ref{2.8}). However, to simplify the computations, we took those values
of $\nabla u\left( x,\overline{k}\right) \mid _{\partial \Omega }$ in our
numerical studies, which were computed when solving the forward problem (\ref%
{8.1}).

\textbf{Remark 8.1}. We have observed in our computations that the solution
of the problem (\ref{8.300}) provides an important piece of information.
Indeed, disks surrounding points of the local maxima of $\left\vert \partial
_{x_{3}}V_{0}\left( x\right) \right\vert $ at $x\in P_{-R+\varepsilon }$ for
a small $\varepsilon >0$ accurately indicate $x_{1},x_{2}$ positions of
inclusions, which we are trying to image, see the text below as well as
Figures 1(f)-3(f).

Now, to update tail functions, we need to follow step 2(b)iii of Algorithm
5.1. More precisely, we need to solve equation (\ref{8.8}) and then use
formula (\ref{4.250}) for $\nabla V_{n,i}.$ However, to speed up
computations, we have decided to use the data $g\left( x,\overline{k}\right) 
$ for this. More precisely, we assume that our inhomogeneities are located
so far from the part $\partial \Omega \diagdown P_{-A}$ of the boundary $%
\partial \Omega $ that their presence provides only very small impact on
this part of the boundary, as compared to their impact on $P_{-A}$. Hence,
we approximately impose the same boundary conditions on $\partial \Omega
\diagdown P_{-A}$ as ones in second and third lines of (\ref{8.1}). Thus,
find the function $u_{n,i}\left( x,\overline{k}\right) $ as the FEM solution
of the following boundary value problem:%
\begin{equation}
\begin{array}{c}
\Delta u_{n,i}+k^{2}c_{n,i}\left( x\right) u_{n,i}=0\text{ in }\Omega , \\ 
u_{n,i}=g\left( x,\overline{k}\right) ,x\in P_{-A}, \\ 
\partial _{n}u_{n,i}+i\overline{k}u_{n,i}=0,x\in \partial _{2}\Omega , \\ 
\partial _{n}u_{n,i}=0,x\in \partial _{3}\Omega .%
\end{array}
\label{8.301}
\end{equation}%
Next, we use formula (\ref{4.250}) to calculate $\nabla V_{n,i}.$ The
question of the well-posedness of problem (\ref{8.301}) is outside of the
scope of this publication. In our computations we did not observe any signs
of the ill-posedness.

\subsubsection{Computations of $c_{n,i}\left( x\right) $}

\label{sec:8.3.2}

It follows from (\ref{4.25}) and (\ref{4.26}) that the function $\beta
_{n,i}\left( x\right) =c_{n,i}\left( x\right) -1$ should be calculated via
applying finite differences to the function $\nabla v_{n,i}\left( x\right) $
given by (\ref{4.24}). The software FreeFem++ automatically interpolates any
function, defined by finite elements to the rectangular grid and we use this
grid to arrange finite differences. Our grid step size is 0.2. As it was
pointed out in Remark 8.1, we have observed in our computations that disks
surrounding the local maxima of the function $\left\vert \partial
_{x_{3}}V_{0}\left( x\right) \right\vert $ at $x\in P_{-A+\varepsilon }$ for
a small $\varepsilon >0$ provide accurate $x_{1},x_{2}$ coordinates of
positions of abnormalities, which we image. Let $x_{1,0}$ and $x_{2,0}$ be $%
x_{1},x_{2}$ coordinates of that point of a local maximum. Then we consider
the cylinder 
\begin{equation}
Cr=\left\{ \left( x_{1},x_{2},x_{3}\right) :\left( x_{1}-x_{1,0}\right)
^{2}+\left( x_{2}-x_{2,0}\right) ^{2}<r^{2},x_{3}\in \left( -A,A\right)
\right\} ,  \label{8.302}
\end{equation}%
where the radius $r=0.3.$ Let $\widetilde{\beta }_{n,i}\left( x\right) $ be
the function computed by the right hand side of (\ref{4.25}). This function
might attain complex or negative values at some points. But we need $\beta
\left( x\right) \geq 0,$ see (\ref{2.2}). Nevertheless, we observed that the
maximal value of the real part of $\widetilde{\beta }_{n,i}\left( x\right) $
in each cylinder (\ref{8.302}) is always positive. Hence, assume that we
have $l$ cylinders $\left\{ \left( Cr\right) _{j}\right\} _{j=1}^{l}$ and
let $\left( \overline{Cr}\right) _{j_{1}}\cap \left( \overline{Cr}\right)
_{j_{2}}$ $=\varnothing $ if $j_{1}\neq j_{2}.$ Then we use the following
truncation to get the function $\beta _{n,i}\left( x\right) :$%
\begin{equation*}
\widehat{\beta }_{n,i}\left( x\right) =\left\{ 
\begin{array}{c}
\func{Re}\widetilde{\beta }_{n,i}\left( x\right) ,\text{ if }\func{Re}%
\widetilde{\beta }_{n,i}\left( x\right) >0.35\max_{\left( \overline{Cr}%
\right) _{j}}\left( \widetilde{\beta }_{n,i}\left( x\right) \right) \text{
and }x\in \left( \overline{Cr}\right) _{j}, \\ 
\text{ for a }j\in \left[ 1,l\right] , \\ 
0,\text{ if either }\func{Re}\widetilde{\beta }_{n,i}\left( x\right) \leq
0.35\max_{\left( \overline{Cr}\right) _{j}}\left( \widetilde{\beta }%
_{n,i}\left( x\right) \right) \text{ or }x\notin \cup _{j=1}^{l}\left( 
\overline{Cr}\right) _{j}\text{ }.%
\end{array}%
\right.
\end{equation*}%
In order to refine images, we have averaged computed functions $\widehat{%
\beta }_{n,i}\left( x\right) $ at each grid point of that rectangular grid.
For each such point we have used nineteen (19) points for averaging: one
point is that grid point and six (6) neighboring points of that rectangular
grid in each of three directions $x_{1},x_{2},x_{3}$. This way we have
obtained the function $\beta _{n,i}\left( x\right) .$ Next, we use (\ref%
{4.26}) to set $c_{n,i}\left( x\right) =\beta _{n,i}\left( x\right) +1.$%
\textbf{\ }

\subsection{Numerical experiments}

\label{sec:8.4}

In this subsection, we present results of our numerical experiments. We
specify domains $\Omega $ and $G$ as 
\begin{equation*}
\Omega =(-2.5,2.5)^{3},G=G=(-3,3)^{3}.
\end{equation*}%
Hence, the part of the boundary $\partial \Omega $ where the backscattering
data $g\left( x,k\right) $ are given, is $P_{-2.5},$ see (\ref{8.0}).
Regardless on the smoothness condition (\ref{2.0}), we reconstruct functions 
$c^{\ast }(x)$ here in the form of step functions. So in each numerical
experiment the support of the function $\beta ^{\ast }\left( x\right)
=c^{\ast }(x)-1$ is in either one or two small inclusion. Thus, 
\begin{equation*}
c^{\ast }(x)=\left\{ 
\begin{array}{c}
3\text{ if }x\text{ is in a small inclusion,} \\ 
1,\text{ otherwise.}%
\end{array}%
\right.
\end{equation*}%
Hence, the inclusion/background contrast is 3 in all cases. We note that
computations usually provide results under lesser restrictive condition than
the theory. In fact, the above mentioned previous results for time dependent
data, including experimental data of \cite{BK1,KFB,KSNF1,TBKF1,TBKF2}, were
also obtained without obeying similar smoothness conditions. Indeed, it is
hard to arrange in experiments such inclusions, which, being embedded in a
medium, would represent, together with that medium, a smooth function. Thus,
this comes back to the point mentioned in the beginning of section 8: about
some discrepancies between the theory and its numerical implementation.

\begin{figure}[h!]
\subfloat[\label{true1a}]{
		\includegraphics[width=\width, height = \height]{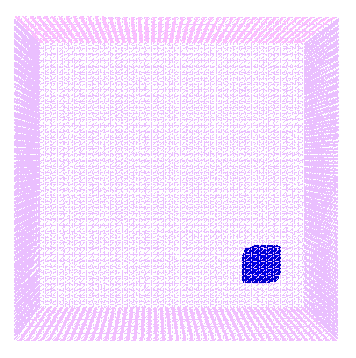} 	} 
\subfloat[\label{true1b}]{
		\includegraphics[width=\width, height = \height]{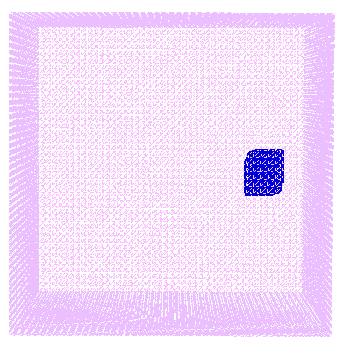} 	} 
\subfloat[]{
		\includegraphics[width=\width, height = \height]{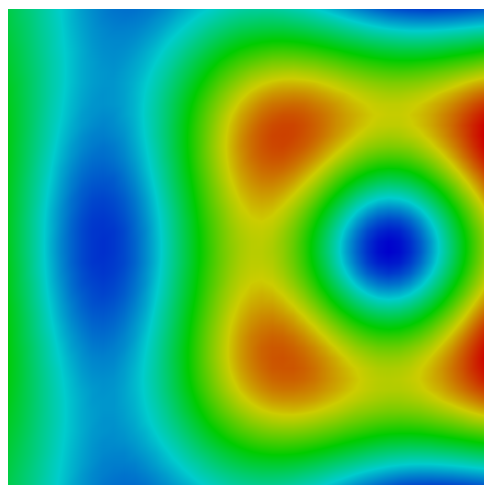} }
\par
\subfloat[]{
		\includegraphics[width=\width, height = \height]{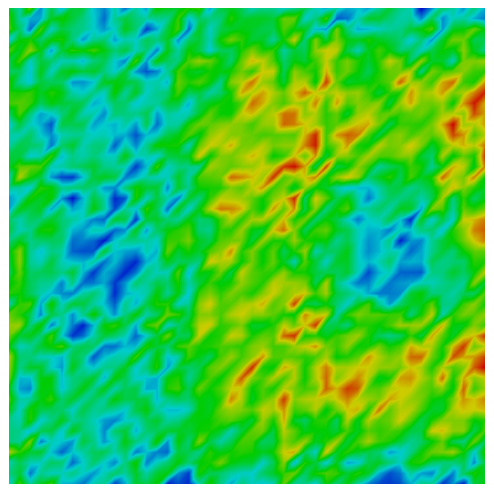} } 
\subfloat[]{
		\includegraphics[width=\width, height = \height]{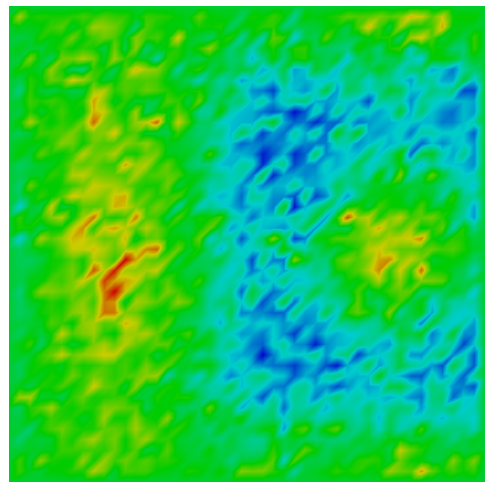} } 
\subfloat[]{
		\includegraphics[width=\width, height = \height]{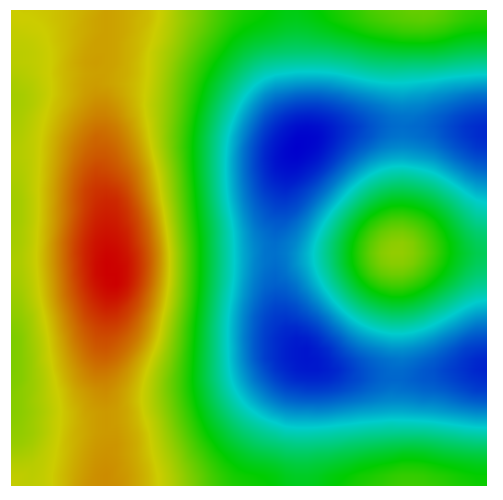} }
\par
\subfloat[]{
		\includegraphics[width=\width, height = \height]{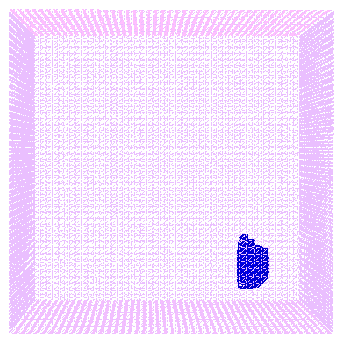} } 
\subfloat[]{
		\includegraphics[width=\width, height = \height]{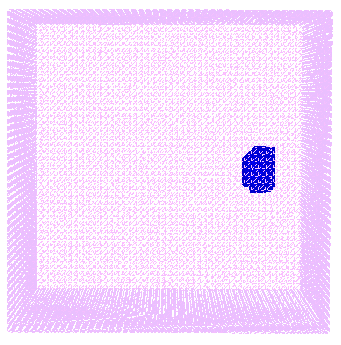} }
\caption{Case 1. $c_{comp}\left( x\right) =c_{7,2}\left( x\right) ,$ $\max
c_{comp}\left( x\right) =2.99.$ The relative error in the maximal value is $%
\left( 3/2.99-1\right) \cdot 100\%=0.33\%.$}
\label{fig 1}
\end{figure}

\begin{figure}[h!]
\subfloat[\label{true2a}]{
		\includegraphics[width=\width, height = \height]{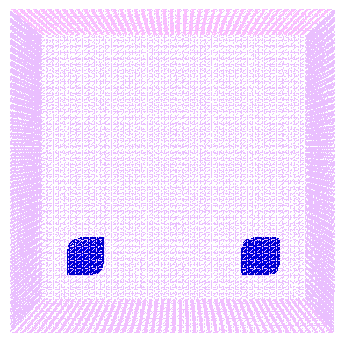} 	} 
\subfloat[\label{true2b}]{
		\includegraphics[width=\width, height = \height]{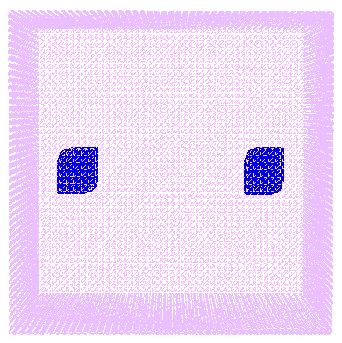} 	} 
\subfloat[]{
		\includegraphics[width=\width, height = \height]{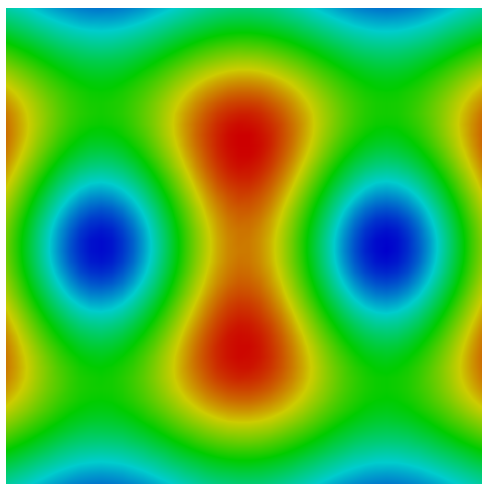} }
\par
\subfloat[]{
		\includegraphics[width=\width, height = \height]{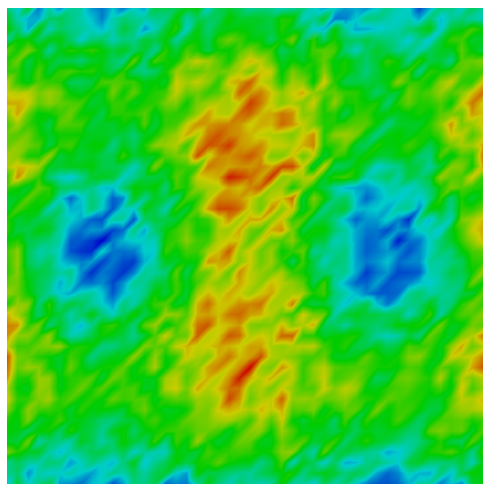} } 
\subfloat[]{
		\includegraphics[width=\width, height = \height]{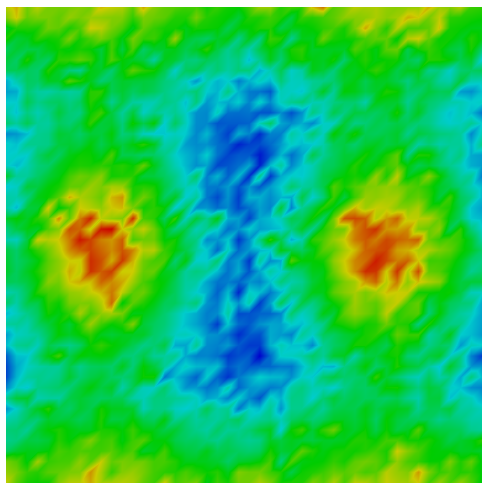} } 
\subfloat[]{
		\includegraphics[width=\width, height = \height]{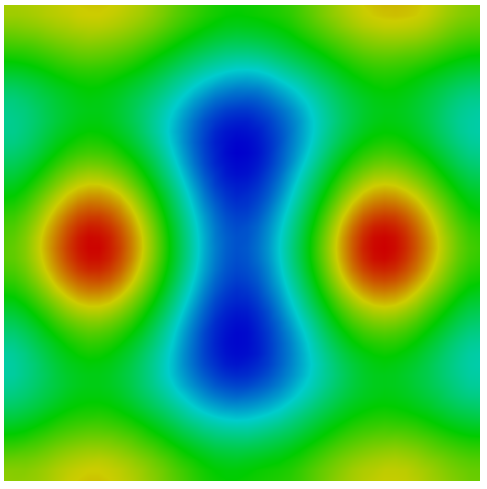} }
\par
\subfloat[]{
		\includegraphics[width=\width, height = \height]{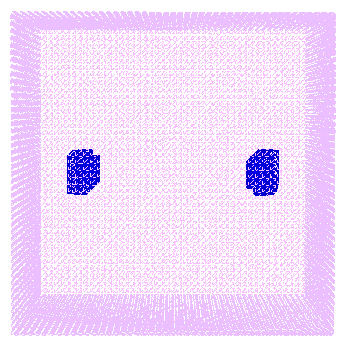} } 
\subfloat[]{
		\includegraphics[width=\width, height = \height]{F2h} }
\caption{Case 2. $c_{comp}\left( x\right) =c_{7,2}\left( x\right) ,$ $\max
c_{comp}\left( x\right) =3.11$. The relative error in the maximal value is $%
3.67\%.$.}
\label{fig 2}
\end{figure}

\begin{figure}[h!]
\subfloat[\label{true3a}]{
		\includegraphics[width=\width, height = \height]{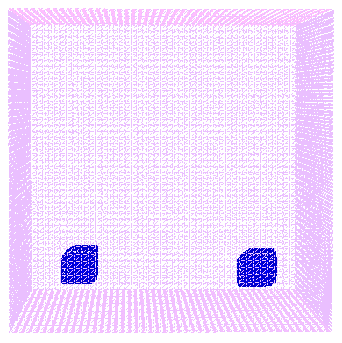} 	} 
\subfloat[\label{true3b}]{
		\includegraphics[width=\width, height = \height]{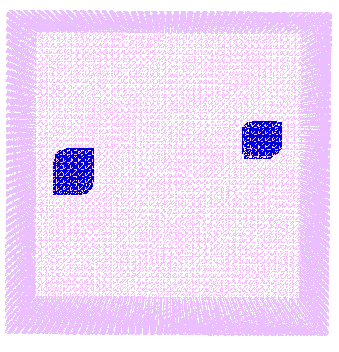} 	} 
\subfloat[]{
		\includegraphics[width=\width, height = \height]{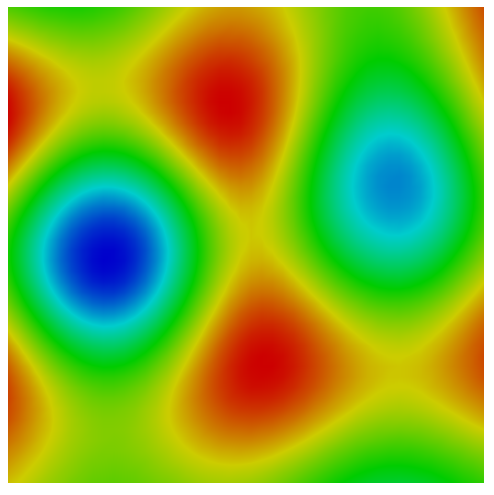} }
\par
\subfloat[]{
		\includegraphics[width=\width, height = \height]{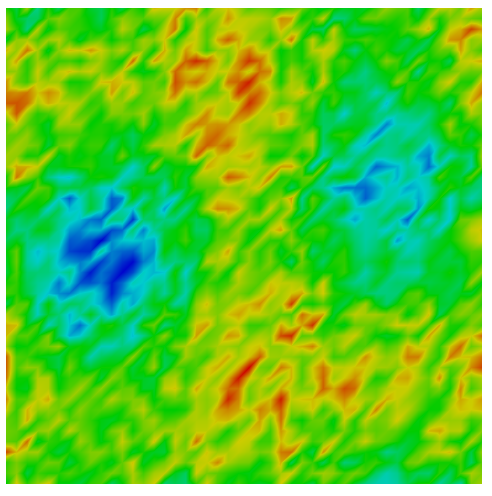} } 
\subfloat[]{
		\includegraphics[width=\width, height = \height]{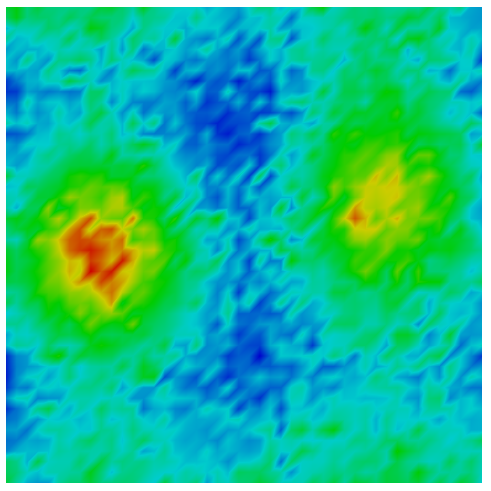} } 
\subfloat[]{
		\includegraphics[width=\width, height = \height]{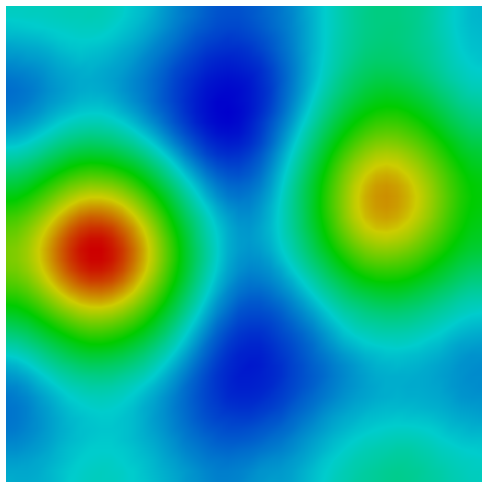} }
\par
\subfloat[]{
		\includegraphics[width=\width, height = \height]{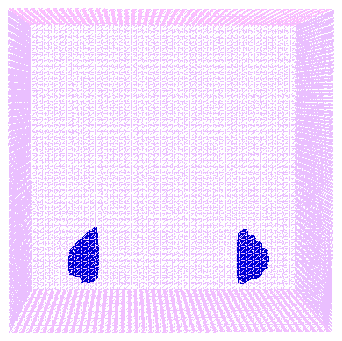} } 
\subfloat[]{
		\includegraphics[width=\width, height = \height]{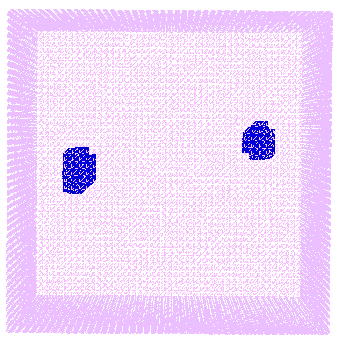} }
\caption{Case 3. $c_{comp}\left( x\right) =c_{7,2}\left( x\right) ,$ $\max
c_{comp}\left( x\right) =3.00$. The relative error in the maximal value is $%
0\%.$}
\label{fig 3}
\end{figure}

In our numerical experiments we test three cases. Inclusions are cubes in
all three. The length of the side of each such cube is 0.5. Our three cases
are:

\begin{enumerate}
\item \textbf{Case 1}. One inclusion. It is centered at $(0,1.5,-1.5)$. See
Figures \ref{true1a}, \ref{true1b}.

\item \textbf{Case 2}. Two inclusions. They are symmetric with respect to
the plane $\{x_{2}=0\}$. The centers of these cubes are at $(0,-1.5,-1.5)$
and at $(0,1.5,-1.5)$. See Figures \ref{true2a}, \ref{true2b}.

\item \textbf{Case 3}. Two inclusions located non-symmetrically with respect
to each of coordinate planes. Their centers are at $(0,-1.5,-1.5)$ and at $%
(1,1.5,-1.5)$. See Figures \ref{true3a}, \ref{true3b}.
\end{enumerate}

We have chosen the $k-$interval as $k\in \lbrack \underline{k},\overline{k}%
]=[1,2].$ Even though our above analysis is valid only for sufficiently
large values of $\underline{k},\overline{k},$ actually it is not clear in
real computations which specific values of these parameters are indeed
sufficiently large. The main reason of our choice of the interval $[%
\underline{k},\overline{k}]$ is that the solution of the problem (\ref{2.3}%
)-(\ref{2.5}) is highly oscillatory for large values of $k$, due to the
presence of the function $u_{0}\left( x,k\right) =\exp \left(
-ikx_{3}\right) .$ So, it takes a lot of computational effort to work with
this solution then. The latter, however, is not the main topic of this
paper, although we will likely study this topic with more details in the
future.

In each of the above three cases, we have chosen $h=0.1$ for the step size
with respect to $k\in \lbrack 1,2]$. Hence, $n=0,\cdots ,N=9.$ In each case,
we generate the data $g\left( x,k_{n}\right) ,n=0,...,9$ for $x\in P_{-2.5}$
via solving the forward problem (\ref{8.1}) for the function $u(x,k_{n})$
and then set $u(x,k_{n})\mid _{P_{-2.5}}=g\left( x,k_{n}\right) .$ We also
add random noise to the data $g\left( x,k_{n}\right) $. The level of this
noise is $5\%$. More precisely, we introduce the noise as 
\begin{equation}
g_{noisy}(x,k_{n})=g(x,k_{n})(1+0.05(\sigma _{1,n}(x)+i\sigma _{2,n}(x))).
\label{8.303}
\end{equation}%
Here $x$ is any vertex of our finite element grid and where $\sigma
_{1,n}(x) $ and $\sigma _{2,n}(x)$ are random numbers in $\left[ -1,1\right] 
$, which are generated by FreeFem++. By (\ref{4.7}) we need to approximate
the derivative $\partial _{k}g_{noisy}(x,k_{n}).$ The differentiation of a
noisy function is an ill-posed problem. So, in our specific case we use a
simple procedure to for the differentiation, 
\begin{equation}
\partial _{k}g_{noisy}\left( x,k_{n}\right) =\frac{g_{noisy}\left(
x,k_{n}\right) -g_{noisy}\left( x,k_{n}-h\right) }{h}.  \label{8.304}
\end{equation}%
We have not observed any instability in this case. This is probably because
the grid step size can sometimes be considered as a regularization parameter
of the differentiation procedure \cite{BK1} and probably our step size $%
h=0.1 $ was suitable for our specific case. However, it is outside of the
scope of this paper to study this question in detail.

Hence, it follows from (\ref{4.00}) and (\ref{8.304}) that we can use only
eight (8) values of $k$: $k_{0}=\overline{k}=2,k_{1}=1.9,...,k_{8}=1.1.$ As
to the number of iterations, our computational experience has shown to us
that the optimal choice was $m=2$, $\overline{N}=7.$ Hence, our computed
functions $c\left( x\right) $ are $c_{comp}\left( x\right) :=c_{7,2}\left(
x\right) $ in all three cases.

Figures \ref{fig 1}, \ref{fig 2} and \ref{fig 3} display our numerical
results for above cases 1, 2 and 3 respectively. In each of these figures we
present:

\begin{enumerate}
\item[(a)] The front view of $\Omega $ for the true model. The data $%
g(x,k_{n})$ are given at the bottom side of $\Omega .$

\item[(b)] The bottom view of $\Omega $ for the true model.

\item[(c)] The absolute value $\left\vert g(x,\overline{k})\right\vert $ of
noiseless data on the measurement plane, which is the bottom side of the
cube $\Omega ,$ i.e. for $x\in P_{-2.5}$.

\item[(d)] The absolute value $\left\vert g_{noisy}(x,\overline{k}%
)\right\vert $ of the noisy data for $x\in P_{-2.5}.$ The disk-like areas
surrounding local minimizers of absolute values of noiseless and noisy data
in Figures \ref{fig 1} (c), (d) and \ref{fig 2} (c), (d) accurately provide $%
x_{1},x_{2}$ positions of true inclusions. However, this cannot be seen
clearly on Figure \ref{fig 3} (d).

\item[(e)] The calculated $|\partial _{x_{3}}V_{0}|$ on the measuring plane $%
P_{-2.5}.$

\item[(f)] The calculated $|\partial _{x_{3}}V_{0}|$ on the plane $P_{-2.4}$%
. We observe that in all Figures 1 (f)-3 (f), the neighborhoods of local
maximizers of $|\partial _{x_{3}}V_{0}|$ accurately provide the $x_{1},x_{2}$
positions of the true inclusions.

\item[(g)] The front view of $\Omega $ for the computed target coefficient $%
c_{comp}(x)=c_{7,2}\left( x\right) $.

\item[(h)] The bottom view of $\Omega $ for the computed target coefficient $%
c_{comp}(x)=c_{7,2}\left( x\right) $.
\end{enumerate}

We observe that the relative errors in maximal values of the function $%
c_{comp}(x)$ are very small. For comparison, we also mention here results
for experimental time dependent data, which were obtained by the globally
convergent numerical method of the first type \cite{BK1,KFB}. Experimental
data of \cite{BK1,KFB} are much noisier of course than our case of (\ref%
{8.303}), see Figures 5.2-5.4 in \cite{BK1} and Figures 3-5 in \cite{KFB}.
Still, Table 5.5 of \cite{BK1} and Table 6 of \cite{KFB} show that the
relative errors in maximal values of the function $c_{comp}(x)$ were varying
between 0.56\% and 2.8\% in four (4) out of five (5) available cases, and
that error was 7.8\% in the fifth case.

Another interesting observation here is that shapes of inclusions are imaged
rather accurately, at least their convex hulls. On the other hand, in the
case of the globally convergent method of \cite{BK1,KFB,KSNF1,TBKF1,TBKF2},
only locations of abnormalities and maximal values of the function $%
c_{comp}(x)$ in them were accurately imaged. So, to image shapes, a locally
convergent Adaptive Finite Element Method was applied on the second stage of
the imaging procedure, see, e.g. Chapters 4 and 5 in \cite{BK1}. We believe
that the better quality of images of shapes here is probably due to a better
sensitivity of the frequency dependent data, as compared with the
sensitivity of the Laplace transformed data.

\section{Summary}

\label{sec:9}

The globally convergent numerical method of the first type, which was
previously developed in \cite{BK1,KFB,KSNF1,TBKF1,TBKF2}, is extended to the
case of the frequency dependent data. The algorithm is developed and its
global convergence is proved. Our method is numerically implemented and
tested for the case of backscattering noisy data. Computational results
demonstrate quite a good accuracy of this technique in imaging of locations
of inclusions, maximal values of the target coefficient $c\left( x\right) $
in them and their shapes.

\begin{center}
\textbf{Acknowledgments}
\end{center}

This work was supported by US Army Research Laboratory and US Army Research
Office grant W911NF-15-1-0233 and by the Office of Naval Research grant
N00014-15-1-2330.

\end{document}